\documentclass[UTF-8,reqno]{amsart}
\usepackage{enumerate, bbm, framed}
\usepackage{pgfplots}
\usetikzlibrary{plotmarks, pgfplots.groupplots}
\usepackage{amssymb,url,color, booktabs}
\usepackage{tikz}
\usepackage{pgfplots}
\pgfplotsset{compat=1.17}
\usetikzlibrary{intersections}
\usepgfplotslibrary{fillbetween}

\usepackage{caption}
\usepackage{float}

\usepackage{mathrsfs}
\usepackage{dutchcal}
\usepackage{threeparttable}

\usepackage{colortbl}
\usepackage{color}
\usepackage[colorlinks=true]{hyperref}
\hypersetup{
    linkcolor=blue,          
    citecolor=red,        
    filecolor=blue,      
    urlcolor=cyan
}

\usepackage{graphicx}
\usepackage{subcaption}
\usepackage{caption}
\usepackage{lipsum}

\usepackage{pgfplots}
\pgfplotsset{compat=newest}
\usepackage{color}
\usepackage{ulem}

 \usepackage{scalerel} 

\definecolor{MyDarkBlue}{cmyk}{0.8,0.3,0.8,0.4}
\definecolor{yellow}{rgb}{0.99,0.99,0.70}
\definecolor{white}{rgb}{1.0,1.0,1.0}
\definecolor{black}{rgb}{0.00,0.00,0.00}

\numberwithin{equation}{section}

\newcommand{\be}{\begin{eqnarray}}
\newcommand{\ee}{\end{eqnarray}}
\newcommand{\ce}{\begin{eqnarray*}}
\newcommand{\de}{\end{eqnarray*}}
\newtheorem{theorem}{Theorem}[section]
\newtheorem{lemma}[theorem]{Lemma}
\newtheorem{remark}[theorem]{Remark}
\newtheorem{definition}[theorem]{Definition}
\newtheorem{proposition}[theorem]{Proposition}
\newtheorem{Examples}[theorem]{Example}
\newtheorem{corollary}[theorem]{Corollary}

\usepackage[nobysame]{amsrefs}
\BibSpec{article}{%
+{}{\PrintAuthors} {author}
+{,}{ \textrm} {title}
+{.}{ \textit} {journal}
+{,}{ \textbf} {volume}
+{}{ \parenthesize} {date}
+{,}{ } {pages}
+{.}{ arXiv:} {eprint}
+{.}{} {transition}
}
\BibSpec{book}{%
+{}{\PrintAuthors} {author}
+{,}{ \textit} {title}
+{.}{ \textrm} {series} 
+{,}{ Vol.} {volume} 
+{.}{ } {publisher}
+{,}{ } {date}
+{.}{} {transition}
}

\def\eps{\varepsilon}

\def\e{\mathrm{e}}

\def\p{\partial}

\def\[{{\Big[}}
\def\]{{\Big]}}
\def\<{{\langle}}
\def\>{{\rangle}}
\def\({{\Big(}}
\def\){{\Big)}}

\def\bx{{\mathbf{x}}}
\def\tr{\mathrm {tr}}

\def\dif{{\mathord{{\rm d}}}}

\def\min{{\mathord{{\rm min}}}}

\def\no{\nonumber}
\def\={&\!\!=\!\!&}

\def\cD{{\mathcal D}}

\def\cI{{\mathcal I}}

\def\cN{{\mathcal N}}

\def\cP{{\mathcal P}}

\def\cS{{\mathcal S}}

\def\cW{{\mathcal W}}
\def\cX{{\mathcal X}}

\def\mE{{\mathbb E}}

\def\mI{{\mathbb I}}

\def\mN{{\mathbb N}}

\def\mP{{\mathbb P}}

\def\mR{{\mathbb R}}

\def\1{{\mathbf{1}}}

\def\sI{{\mathscr I}}

\def\geq{\geqslant}
\def\leq{\leqslant}

\def\div{\mathord{{\rm div}}}

\def\eps{\varepsilon}

\def\e{\mathrm{e}}

\def\p{\partial}

\def\[{{\Big[}}
\def\]{{\Big]}}
\def\<{{\langle}}
\def\>{{\rangle}}
\def\({{\Big(}}
\def\){{\Big)}}

\def\bx{{\mathbf{x}}}
\def\tr{\mathrm {tr}}

\def\dif{{\mathord{{\rm d}}}}

\def\min{{\mathord{{\rm min}}}}

\def\no{\nonumber}
\def\={&\!\!=\!\!&}
\def\bt{\begin{theorem}}
\def\et{\end{theorem}}
\def\bl{\begin{lemma}}
\def\el{\end{lemma}}
\def\br{\begin{remark}}
\def\er{\end{remark}}
\def\bx{\begin{Examples}}
\def\ex{\end{Examples}}
\def\bd{\begin{definition}}
\def\ed{\end{definition}}
\def\bp{\begin{proposition}}
\def\ep{\end{proposition}}
\def\bc{\begin{corollary}}
\def\ec{\end{corollary}}

\def\geq{\geqslant}
\def\leq{\leqslant}

\def\div{\mathord{{\rm div}}}

\def\<{\langle} \def\>{\rangle}

\def\bpf{\begin{proof}}
\def\epf{\end{proof}}

\allowdisplaybreaks

\begin{document}

\title{New algorithms for sampling and diffusion models}
\author{Xicheng Zhang}

\address{Xicheng Zhang:
School of Mathematics and Statistics, Beijing Institute of Technology, Beijing 100081, China\\
Email: xczhang.math@bit.edu.cn
 }

\thanks{{\it Keywords: \rm Sampling, reversed diffusion process, probability measure flow}}

\thanks{
This work is supported by National Key R\&D program of China (No. 2023YFA1010103) and NNSFC grant of China (No. 12131019)  and the DFG through the CRC 1283 ``Taming uncertainty and profiting from randomness and low regularity in analysis, stochastics and their applications''. }


\begin{abstract}
Drawing from the theory of stochastic differential equations, we introduce a novel sampling method for known distributions and a new algorithm for diffusion generative models with unknown distributions. Our approach is inspired by the concept of the reverse diffusion process, widely adopted in diffusion generative models. Additionally, we derive the explicit convergence rate based on the smooth ODE flow. 
For diffusion generative models and sampling, we establish a {\it dimension-free} particle approximation convergence result.
Numerical experiments demonstrate the effectiveness of our method. Notably, unlike the traditional Langevin method, our sampling method does not require any regularity assumptions about the density function of the target distribution. Furthermore, we also  apply our method to optimization problems. 

\end{abstract}

\maketitle
\setcounter{tocdepth}{2}

\tableofcontents

\section{Introduction}

\subsection{Background}

In statistics, machine learning, and data science, two fundamental challenges are commonly encountered: drawing samples from a known distribution and learning an unknown distribution from observed data to generate new samples.

For the former problem, when we have a well-defined probability distribution (e.g., Gaussian, binomial, exponential) with known parameters, we can generate random samples that follow this distribution. This process is essential for tasks like Monte Carlo simulations, hypothesis testing, and generating synthetic data for modeling purposes. Nowadays, various algorithms exist for sampling from known distributions, depending on the nature of the distribution. For example, the inverse transform method, rejection sampling, and Markov chain Monte Carlo (MCMC) methods such as the Metropolis-Hastings algorithm and Gibbs sampling are commonly used techniques for generating samples from different types of distributions (see \cite{L10, RK16}).

For the latter problem, the goal is to infer the underlying probability distribution that generated a given dataset without prior knowledge of the distribution's form. This process involves using statistical techniques and algorithms to estimate the parameters or structure of the distribution based on the observed data. Methods for learning unknown distributions include parametric and non-parametric approaches. Parametric methods assume a specific functional form for the distribution (e.g., Gaussian, Poisson) and estimate the parameters that best fit the data. Non-parametric methods, on the other hand, make fewer assumptions about the distribution's shape and instead focus on estimating the distribution directly from the data, often using techniques like kernel density estimation or histogram methods (see \cite{GCSR04, M12, HTF08, B06}).

In recent years, diffusion generative models have emerged as powerful frameworks for generating high-quality synthetic data, 
particularly in the realm of image generation. These models operate by iteratively transforming a simple noise distribution into 
a complex data distribution, effectively reversing a diffusion process.
Diffusion generative models leverage the principles of probabilistic modeling and stochastic processes to generate data that closely mimics real-world distributions. Introduced by Sohl-Dickstein et al. in 2015 \cite{DWMG15} 
and popularized by subsequent works such as those by Ho et al. in 2020 \cite{HJA20}, 
these models have demonstrated remarkable success in various generative tasks, particularly in generating high-quality images.

In these contexts, two notable probabilistic generative models have been developed: Score Matching with Langevin Dynamics (SMLD), introduced by Song and Ermon in 2019 \cite{SE19}, and Denoising Diffusion Probabilistic Modeling (DDPM), introduced by Ho et al. in 2020 \cite{HJA20}. 
Both models provide robust frameworks for generative modeling.
The SMLD method estimates the score function (i.e., the gradient of the log probability density with respect to the data) at each noise scale and uses Langevin dynamics to sample from a sequence of decreasing noise scales during generation. On the other hand, the DDPM method aims to capture the data distribution by progressively transforming a noise distribution into the target data distribution through an iterative diffusion process.

Since their introduction, numerous studies have explored and expanded upon diffusion generative models. Notably, Song et al. in 2021 \cite{SD21}
proposed a general score-based generative model using stochastic differential equations (SDEs) to unify the methods of SMLD and DDPM. 
The crucial idea is to learn the score function $(t,x)\mapsto\nabla_x p_t(x)$ through a forward SDE by slowly injecting noise into the data. 
This approach allows solving a reverse-time SDE (see \cite{And82} and \cite{HP86}) to generate high-quality samples starting from a normal distribution.
See below for an explanation of the basic ideas (adapted from \cite[Figure 1]{SD21}).
\begin{figure}[ht]
    \centering
 \includegraphics[width=3.5in, height=1.5in]{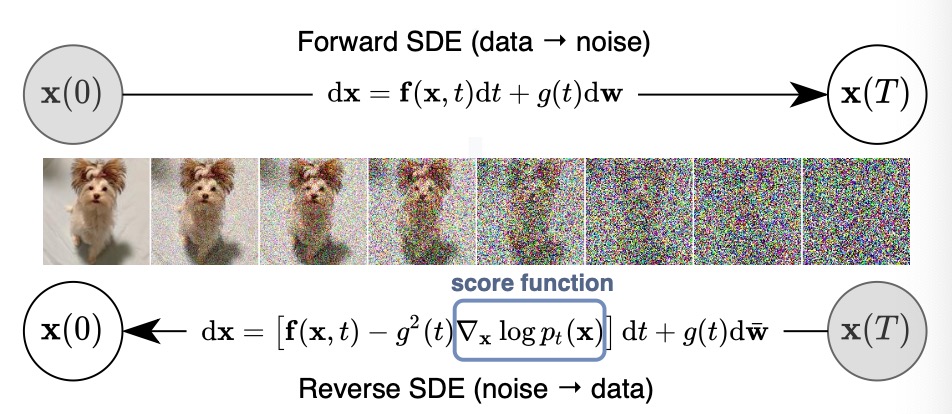}
    \caption{Explanation of diffusion models}
    \label{fig:sd21}
\end{figure}

To learn the score function using neural networks, well-known methods typically employ black-box approaches that involve constructing suitable loss functions. Recent advancements have also introduced new algorithms for sampling and improved diffusion models, further enhancing the quality and efficiency of generative modeling (cf. \cite{LZ22} and the survey paper \cite{YZS23}). Additionally, several works focus on the convergence analysis of diffusion generative models (see \cite{LLT23, CLL23, TL23, CLTX23}).

\subsection{Main results}

Inspired by the reversed diffusion process, we propose a more explicit probability flow ordinary differential equation (ODE) as in \cite{SD21} to generate samples. The key observation is that through the probability flow of the ODE, we can establish a connection between any probability distribution and the standard normal distribution. Since the coefficient of our ODE depends explicitly on the distribution, our method can be considered a white-box approach, in contrast to methods that learn the score function based on data.

More precisely, let $\mu_0$ be any given distribution in $\mathbb{R}^d$ and $\eta \sim \mu_0$. 
To generate a sample, it suffices to solve the following probability flow ODE:
\begin{align}\label{ODE}
\dot{Y}_t = \frac{\mathbb{E}\left[(\eta - Y_t) \exp\left(-\frac12\Big\|\frac{Y_t - t\eta}{1 - t}\Big\|_2^2 \right)\right]}
{(1 - t)\mathbb{E}\left[\exp\left(-\frac12\Big\|\frac{Y_t - t\eta}{1 - t}\Big\|_2^2 \right)\right]} =: b^{\mu_0}_t(Y_t), \quad Y_0 \sim N(0, \mathbb{I}_d),
\end{align}
where the initial value $Y_0$ is an independent $d$-dimensional standard normal distribution, for $q \in [1, \infty]$, $\|\cdot\|_q$ stands for the usual $\ell^q$-norm in $\mathbb{R}^d$, and the expectation is taken only with respect to $\eta$. In particular, we show that the law of $Y_t$ weakly converges to $\mu_0$ as $t \uparrow 1$.

The following theorem is a combination of Lemma \ref{Le27} and Theorem \ref{Th24} below.
\bt
Suppose that $\mu_0\in\cP(\mR^d)$ has compact support and $\eta\sim\mu_0$. Then for each $t\in[0,1)$,
$y\mapsto b^{\mu_0}_t(y)$ is $C^\infty$-smooth and for all $x,y\in\mR^d$,
\begin{align}\label{Mon1}
\<x-y, b^{\mu_0}_t(x)-b^{\mu_0}_t(y)\>\leq \left[\frac{t\sup_\omega\|\eta(\omega)\|^2_2}{(1-t)^2}-1\right]\frac{\|x-y\|^2_2}{1-t},
\end{align}
and for each $Y_0\sim N(0,\mI_d)$ independent of $\eta$,
there is a unique solution $Y_t$ to ODE \eqref{ODE} on $[0,1)$. 
Moreover, let $\mu_t$ be the law of $Y_t$ and
$\cW_2$ be the usual Wasserstein metric, then
\begin{align}\label{AS2}
\cW_{2}(\mu_t,\mu_0)\leq \left(\sup_{\omega}\|\eta(\omega)\|_2+\sqrt{d}\right)(1-t),
\end{align}
and for any $M>0$ and $t\in[0,1)$,
\begin{align}\label{AS1}
\mP\left(\|Y_t\|_\infty\geq t\sup_{\omega}\|\eta(\omega)\|_\infty+\sqrt 2 (1-t) M\right)\leq 1-\left(1-\frac{2}{\sqrt\pi}\int_{M}^{\infty}\e^{-x^2}\dif x\right)^d.
\end{align}
\et
\br
From \eqref{AS1}, one sees that as $t$ approaches $1$, with probability 1,
$$
\varlimsup_{t \to 1} \|Y_t\|_\infty \leq \sup_{\omega} \|\eta(\omega)\|_\infty.
$$
An open question here is to show that $Y_t$ almost surely converges to a random variable $\eta \sim \mu_0$ as $t \to 1$, rather than merely exhibiting weak convergence as in \eqref{AS2}. Numerical experiments indicate that when $\mu_0(\dif x) = \sum_{i=1}^N \delta_{a_i}(\dif x)$ is a discrete distribution, for any sample $Y_0 \sim N(0, \mathbb{I}_d)$, $Y_t$ always converges to some point $a_i$ as $t \to 1$. See Section 5 below.
\er

In practical applications, we usually need to consider the particle approximation of ODE \eqref{ODE}:
\begin{align}\label{ODE1}
\dot{Y}^N_t = \frac{\sum_{j=1}^N\left[(\eta_j - Y^N_t) \exp\left(-\frac12\Big\|\frac{Y^N_t - t\eta_j}{1 - t}\Big\|_2^2 \right)\right]}
{(1 - t)\sum_{j=1}^N\left[\exp\left(-\frac12\Big\|\frac{Y^N_t - t\eta_j}{1 - t}\Big\|_2^2 \right)\right]}, \quad Y^N_0=Y_0 \sim N(0, \mathbb{I}_d),
\end{align}
where $\{\eta_j, j=1,\cdots, N\}$ is a sequence of i.i.d. random variables with common distribution $\mu_0$ and it is also independent of $Y_0$. One may ask whether $Y^N_t $ converges to $Y_t$. In this aspect,
we have the following dimension-free  convergence result, 
which has independent interest (see Theorem \ref{Th51}). 
\bt\label{Th31}
Suppose $K:=\sup_{\omega}\|\eta(\omega)\|_2<\infty$. Then it holds that for all $t\in[0,1)$,
$$
\sup_{N\in\mN}\left(\sqrt N\mE\|Y_t-Y^N_t\|_2\right)\leq 2Kt\e^{(Kt/(1-t))^2}.
$$
\et
\br
It appears that the aforementioned result provides the initial quantitative estimate for approximating generative models of diffusions. Although the exponential component grows rapidly as $ t \to 1 $, in practical numerical experiments, 
we typically set $ t < 1 $, such as $ t = 0.9 $. Additionally, through scaling techniques, we can always choose a small $ K $, for instance, $ K = 0.1 $, ensuring that the exponential part remains less than 1.
\er

Now we explain how to sample from a given distribution $\mu_0(\dif x) \propto f(x) \dif x$ using the ODE \eqref{ODE}, where
$$
\text{$f: \mathbb{R}^d \to [0, \infty)$ is bounded by $\Lambda$ and has support in $\{x: \|x\|_2\leq K\}$.}
$$
In this special case, let $\xi$ be the uniformly distributed random variable in 
$\{x: \|x\|_2\leq K\}$, and $\xi'\sim N(0,\mI_d)$. 
Then we can write for $t \in (0, 1)$ 
$$
b^{\mu_0}_t(y) = \frac{\mathbb{E}\left[(\xi - y) \exp\left(-\frac{1}{2} 
\left\|\frac{y - t\xi}{1 - t}\right\|_2^2 \right) f(\xi)\right]}
{(1 - t)\mathbb{E}\left[\exp\left(-\frac{1}{2} \left\|\frac{y - t\xi}{1 - t}\right\|_2^2 \right) f(\xi)\right]},
$$
and by the integration by parts, an alternative expression is
$$
b^{\mu_0}_t(y)=\frac{\mE[(y-\xi') f((y-(1-t)\xi')/t)]}{t\mE f((y-(1-t)\xi')/t)}.
$$
Since it easily generates a sequence of i.i.d. uniformly or normally distributed random variables, we can calculate $b^{\mu_0}_t(y)$ using the Monte Carlo method as shown in \eqref{ODE1}. 
Thus, we can solve a similar ODE as in \eqref{ODE1} using the classical Euler discretization method. Alternatively, under \eqref{Mon1}, we can employ the recently developed Poisson's discretization approximation for the ODE \eqref{ODE1} as presented in \cite{Zh23}. Note that the method in \cite{Zh23} does not require any time regularity of $b^{\mu_0}_t(y)$.

Nowadays, the popular sampling method is the Markov chain Monte Carlo (MCMC), whose stationary distribution is the target distribution, based on ergodicity theory. Note that MCMC sampling is also a key component in the popular simulated annealing technique for both discrete and continuous optimization. Theoretically, the MCMC method requires the target density function to have some smoothness to ensure the iterations converge.
Compared with MCMC, our new method does not make any regularity assumptions about the density function,
and it is easy to be realized. Numerical experiments exhibit robust sampling ability. 

\subsection{Related works}
Since the seminal works of \cite{HJA20} and \cite{SD21}, numerous papers have explored diffusion generative models (for an overview, see the recent survey \cite{YZS23} and \cite{CMFW24}). Notably, Albergo et al. \cite{ABE23} introduced a general framework for stochastic interpolants that unifies various diffusion models. Specifically, let $\rho_0$ and $\rho_1$ be two density functions in $\mathbb{R}^d$, and let $(x_0,x_1)$ be a coupling of $(\rho_0,\rho_1)$, represented as $\mathbb{R}^{2d}$-valued random variables with marginal distributions $\rho_0$ and $\rho_1$. Define $z$ as a normally distributed variable. The stochastic interpolants between $x_0$ and $x_1$ are defined by
$$
x_t := I(t,x_0,x_1) + \gamma(t) z, \quad t \in [0,1],
$$
where $I \in C^2([0,1]; (C^2(\mathbb{R}^d \times \mathbb{R}^d))^d)$ satisfies $I(0,x_0,x_1) = x_0$, $I(1,x_0,x_1) = x_1$, and
$$
|\partial_t I(t,x_0,x_1)| \leq C |x_0 - x_1|,
$$
while $\gamma: [0,1] \to [0,\infty)$ satisfies $\gamma(0) = \gamma(1) = 0$, $\gamma(t) > 0$ for $t \in (0,1)$, and $\gamma^2 \in C^2([0,1])$. Under certain regularity assumptions on $\rho_0$ and $\rho_1$, the authors in \cite[Theorem 1]{ABE23} demonstrated that the density $\rho_t$ of $x_t$ satisfies the following Fokker-Planck equation:
$$
\partial_t \rho + \text{div}(b \cdot \rho) = 0, \quad 
b(t,x) = \mathbb{E}[\dot{x}_t | x_t = x].
$$
Our equation in Lemma \ref{Le23} below represents a special case of the aforementioned general equation. Unlike \cite[Theorem 1]{ABE23}, we do not impose any regularity assumptions on the target distribution.
Following the general framework of \cite{ABE23}, Gao, Huang, and Jiao \cite{GHJ23} studied the Gaussian interpolant flow and proved several analytic results for the associated ODE flows under various assumptions.


About the sampling from a given distribution, in a similar vein, Huang et al. \cite{HJK21} consider the Schrödinger-Föllmer diffusion (also known as the Föllmer process \cite{F84, ELS20}):
$$
\dif Y_t = \nabla\log P_{1-t}f(Y_t)\dif t + \dif W_t, \quad Y_0 = 0,
$$
where $P_t f(x) := \mathbb{E} [f(x + W_t) \e^{\|x + W_t\|_2^2/2}]$. It was shown that the law of $Y_1$ is $\mu_0$.
Huang et al. \cite{HDMZ24} consider the reverse diffusion process:
$$
\dif Y_t = [Y_t + 2\nabla\log p_{T-t}(Y_t)]\dif t + \sqrt{2} \dif W_t,
$$
where $p_t$ is the density of the forward Ornstein-Uhlenbeck process 
$\dif X_t = -X_t \dif t + \sqrt{2} \dif W_t$ with $X_0 \propto \e^{-f}$.
Another approach is based on stochastic localization \cite{GMGD24}, which is:
$$
\dif Y_t = u_t(Y_t)\dif t + \sigma \dif W_t, \quad Y_0 = 0,
$$
where $u_t(y)=y/t+\sigma^2\nabla\log p_{t}(y)$ and $p_t(y)=(2\pi t\sigma^2)^{-d/2}\int_{\mR^d}\e^{-\|y-tx\|^2_2/(2t\sigma^2)}f(x)\dif x$. It was shown that the law of $Y_t$ converges to $\mu_0$ as $t\to\infty$.
The main difference between our sampling method and the above approaches is that we utilize the probability flow of ODEs rather than SDEs. Importantly, we have explicit convergence rates.

Since this paper primarily focuses on theoretical aspects, and my expertise lies more in theory than in practical algorithms, I will not provide a comparison between our algorithm and various well-known algorithms.

\subsection{Organization}
This paper is organized as follows: In Section 2, we first recall the basic idea of diffusion generative models used in \cite{SD21}. Then, we present our main theoretical result, providing quantitative estimates about the probability flow based on ODE, which is important for devising an algorithm. In Section 3, we demonstrate how to generate samples randomly from a family of high-dimensional discretized distributions.
In Section 4,  we apply our theoretical results to the sampling problem from a known distribution, illustrating the effectiveness of our method through various examples, including discontinuous density functions and multimodal distributions.  In Section 5, we discuss potential applications in optimization problems based on sampling the density function.  

\section{Proofs of main results}
Throughout this paper, we fix the dimension $d\in\mN$ and a probability  measure $\mu_0\in\cP(\mR^d)$.
Let $\mI_d$ be the identity matrix of $d\times d$. For $q\in[1,\infty]$,
let $\|\cdot\|_{q}$ be the usual $\ell^q$-norm in $\mR^d$, i.e.,
\begin{align}\label{Dis}
\|x\|_q:=\left(\sum_{i=1}^d|x_i|^q\right)^{1/q}
\end{align}
with the convention that for $q=\infty$, $\|x\|_\infty:=\max\{|x_1|,\cdots,|x_d|\}$. Notice that $\|\cdot\|_2$ is the usual Euclidean norm.
\subsection{Reverse diffusion processes}\label{Sec1}
We  fix two continuous functions
$$
f_t, g_t: (0,1)\to [0,\infty),
$$ 
which will be specified below. Now we consider the following linear SDE in $\mR^d$:
\begin{align}\label{LSDE1}
\dif X_t=-f_tX_t\dif t+\sqrt{2g_t}\dif W_t,\ \ t\in[0,1),
\end{align}
where $X_0\sim\mu_0$ and $W$ is a $d$-dimensional standard Brownian motion independent with $X_0$. The solution of linear SDE \eqref{LSDE1} is explicitly given by
$$
X_t=\e^{-\int^t_0f_s\dif s}X_0+\int^t_0\e^{-\int^t_sf_r\dif r}\sqrt{2g_s}\dif W_s.
$$
In particular, if we set
$$
h_t:=\e^{-\int^t_0f_s\dif s},
$$
then
\begin{align}\label{AF1}
X_t=h_tX_0+h_t\int^t_0\sqrt{2g_s}/ h_s\dif W_s=:h_t X_0+\xi_t.
\end{align}
Clearly, $\xi_t$ is a $d$-dimensional normal distribution with mean zero and covariance matrix
$$
\mE(\xi^i_t\xi^j_t)=\1_{i=j}a^2_t,\ \ i,j=1,\cdots,d,
$$
where
$$
a^2_t:=2 h^2_t\int^t_0g_s/ h^2_s \dif s,\ \ t\in(0,1).
$$
By the definitions of $h_t, a_t$ and the chain rule, it is easy to see that
\begin{align}\label{FG1}
\boxed{ h'_t=-f_th_t,\ \  a_ta'_t=g_t-f_ta^2_t,}
\end{align}
where the prime stands for the derivative of a function with respect to the time variable.

Now if we assume that $h_t$ is strictly decreasing, $a_t$ is strictly increasing, and
\begin{align}\label{AA1}
\lim_{t\downarrow 0}h_t=1,\ \ \lim_{t\uparrow 1}h_t=0,\ \ \lim_{t\downarrow 0}a_t=0,\ \ \lim_{t\uparrow 1}a_t=1,
\end{align}
then by \eqref{AF1}, it holds that in weak sense,
$$
X_t\stackrel{t\downarrow 0}{\to}X_0\sim\mu_0,\ \ X_t\stackrel{t\uparrow 1}{\to}X_1=\mu_1\sim N(0,\mI_d),
$$
that is, for any bounded continuous function $\varphi:\mR^d\to\mR$,
$$
\lim_{t\to 0}\mE\varphi (X_t)=\int_{\mR^d}\varphi(x)\mu_0(\dif x),\ \ \lim_{t\to 1}\mE\varphi (X_t)=\int_{\mR^d}\varphi(x)\mu_1(\dif x).
$$
In other words, $(X_t)_{t\in(0,1)}$ builds a bridge between the initial distribution $\mu_0$ and the standard normal distribution $\mu_1\sim N(0,\mI_d)$.

\begin{minipage}{\textwidth}
\begin{tikzpicture}[baseline, xscale=0.7, yscale=0.5]
\begin{axis}[
    domain=0:1,
    samples=100,
    grid=both,
    axis lines=middle,
    xlabel=$x$,
    ymax=2
]
\addplot [red, thick] {1 + (sin(deg(2*pi*x)) + sin(deg(4*pi*x)))/2};
\begin{scope}[on layer=axis foreground]
    \node at (axis cs: 0.6, 1.7) {$\mu_0\sim1_{x\in(0,1)}\left[1+\frac{\sin(2\pi x)+\sin(4\pi x)}{2}\right]$};
\end{scope}
\end{axis}
\end{tikzpicture}
\begin{tikzpicture}
\node at (2, 2.5) {\tiny $X_0\to X_t\to X_1$};
\draw [blue, line width=0.5,->] (1,2) -- (3,2);
\draw [red, line width=0.5, <-] (1,1) -- (3,1);

\begin{tikzpicture}
   (0,0) -- (5,0);
\end{tikzpicture}

\node at (2, 0.5) {\tiny $X_0\leftarrow X_t\leftarrow X_1$};
\end{tikzpicture}
\begin{tikzpicture}[baseline, xscale=0.7, yscale=0.5]
\begin{axis}[
    domain=-3:3,
    samples=100,
    grid=both,
    axis lines=middle,
    xlabel=$x$,
    ymax=0.4
]
\addplot [blue, thick] {1/sqrt(2*pi)*exp(-x^2/2)};
\begin{scope}[on layer=axis foreground]
    \node at (axis cs: 2, 0.3) {$\mu_1\sim N(0,1)$};
\end{scope}
\end{axis}
\end{tikzpicture}
\captionof{figure}{}
\end{minipage}
In the diffusion generative models, as $t$ varies from $0$ to $1$, $X_0\to X_t\to X_1$ is called the ``adding noise process''. As $t$ varies from $1$ to $0$,
$X_1\to X_{t}\to X_0$  is called the ``denoising process''.

The following lemma is direct by expression \eqref{AF1}.
\bl\label{Le21}
Let $\mu_t$ be the law of $X_t$. For any $q\in[1,\infty]$ and $p\in[1,\infty)$, we have
$$
\cW_{q,p}(\mu_t,\mu_0)\leq(1-h_t)\left(\int_\Omega\|X_0(\omega)\|^p_q \mP(\dif\omega)\right)^{1/p}+a_t\left(\int_\Omega\|X_1(\omega)\|^p_q \mP(\dif\omega)\right)^{1/p},
$$
where $X_0\sim\mu_0$ and $X_1\sim N(0,\mI_d)$, and the Wasserstein metric $\cW_{q,p}$ is defined by
\begin{align}\label{WW1}
\cW_{q,p}(\mu_t,\mu_0):=\inf_{{\rm Law}(\xi)=\mu_t, {\rm Law}(\eta)=\mu_0}\left(\int_\Omega\|\xi(\omega)-\eta(\omega)\|^p_q\mP(\dif\omega)\right)^{1/p}.
\end{align}
\el
\br
When $1-h_t=a_t$, we have the best convergence rate as $t\to 0$. This is also the reason why we choose $\beta_t=1-\sigma_t$ in {\bf (H)} below.
\er
Next we want to find the evolution equation of denoising process $\bar X_t:=X_{1-t}$.
For $\sigma>0$, let
\begin{align}\label{DH1}
\rho_{\sigma}(x):=(2\pi \sigma^2)^{-d/2}\e^{-\|x\|_2^2/(2\sigma^2)}.
\end{align}
Then for each $t\in(0,1)$, $X_t$ admits a density $p_t(x)$ given by
\begin{align}\label{SD1}
p_t(x)=\mE\rho_{a_t}(x-h_tX_0)=\int_{\mR^d}\rho_{a_t}(x-h_ty)\mu_0(\dif y).
\end{align}
Indeed, it follows by \eqref{AF1} that
$$
\int_{\mR^d}\varphi(x)p_t(x)\dif x=\mE \varphi(X_t)=\mE \varphi(h_tX_0+\xi_t)=\int_{\mR^d}\mE \varphi(h_tX_0+x)\rho_{a_t}(x)\dif x.
$$
On the other hand, for any $\varphi\in C^2_b(\mR^d)$,
by \eqref{LSDE1} and It\^o's formula,  we have
$$
\mE \varphi(X_t)=\mE \varphi(X_0)+\int^t_0\mE\big[g_s\Delta \varphi(X_s)-f_sX_s\cdot\nabla \varphi(X_s)\big]\dif s,
$$
and so,
$$
\int_{\mR^d}\varphi(x) p_t(x)\dif x=\int_{\mR^d} \varphi(x)\mu_0(\dif x)
+\int^t_0\!\!\int_{\mR^d}[g_s\Delta \varphi(x)-f_sx\cdot\nabla \varphi(x)]p_s(x)\dif x\dif s.
$$
From this, one sees that $p_t(x)$ solves the following linear Fokker-Planck equation:
$$
\p_tp_t(x)=g_t\Delta p_t(x)+\div (f_tx\, p_t(x)).
$$
For a function $f:[0,1]\to\mR$, we  write
$$
\bar f_t:=f_{1-t}.
$$
By this  notation and the chain rule, it is easy to see that
\begin{align}
\p_t\bar p_t(x)&=-\bar g_t \Delta \bar p_t(x)-\div (\bar f_tx\, \bar p_t(x))\no\\
&=-\div ((\bar g_t\nabla\log\bar p_t(x)+\bar f_tx)\bar p_t(x))\no\\
&=\epsilon_t\Delta \bar p_t(x)-\div (\bar b_t(x)\bar p_t(x)),\label{PDE1}
\end{align}
where $\epsilon_t: [0,1]\to[0,\infty)$ and
\begin{align}\label{AF2}
\bar b_t(x):=(\bar g_t+\epsilon_t)\boxed{\nabla\log \bar p_t(x)}+\bar f_tx.
\end{align}
In the literature of diffusion models, the function $\nabla \log \bar{p}_t(x)$ is called the score function, which depends on the unknown distribution.
By the superposition principle (for example, see \cite[Theorem 1.5]{RXZ20}), formally, $\bar{p}_t(x)$ will be the density of the following SDE:
$$
\dif Y_t = \sqrt{2\epsilon_t}\dif W_t + \bar{b}_t(Y_t)\dif t, \quad Y_0 \sim N(0, \mathbb{I}_d).
$$
In particular, if $\epsilon_t = \bar{g}_t$, then $(Y_t)_{t \in (0, 1)}$ has the same law as the reverse diffusion process $(\bar{X}_t)_{t \in (0, 1)}$ 
(cf. \cite[Theorem 2.1]{And82}).

\subsection{Probability measure flows}

Motivated by the above discussion, 
in what follows, we fix a terminal time $T$ (finite or infinite) and directly consider a family of probability density functions:
\begin{align}\label{Den1}
\phi_t(x):=\int_{\mR^d}\rho_{\sigma_t}(x-\beta_ty)\mu_0(\dif y),\ \ t\in(0,T),
\end{align}
where $\sigma_t,\beta_t: [0,T)\to[0,\infty)$ are two continuous functions and satisfy that
\begin{enumerate}[{\bf (H)}]
\item $\sigma_t$ is stricly decreasing and for any $t_0\in(0,T)$, $\sigma'_t$ is bounded on $[0,t_0]$, and 
$$
\lim_{t\downarrow 0}\sigma_t=1,\ \ \lim_{t\uparrow T}\sigma_t=0,\ \ \beta_t=1-\sigma_t.
$$
\end{enumerate}
\br
When $T=\infty$, we typically choose $\sigma_t=\e^{-t}$, and thus
$$
\phi_t(x)=\int_{\mR^d}\rho_{\e^{-t}}(x-(1-\e^{-t})y)\mu_0(\dif y),\ \ t>0.
$$ 
When $T=1$,
the typical choice of $\sigma_t$ is
$$
\sigma_t=(1-t)^\alpha\ \ \mbox{or}\ \  \sigma_t=1-t^\alpha,\ \ \alpha>0.
$$
Here are the graphs of $\sigma_t$ and $\beta_t$ for $\alpha=0.5$ and $\alpha=2$.
\vspace{5mm}
\\
\begin{tikzpicture}
\begin{groupplot}[group style={group size=2 by 1, horizontal sep=2cm, vertical sep=2cm}]
    \nextgroupplot[
        xlabel={$\alpha=0.5$},
        width=7cm,
        height=5cm,
        legend style={at={(1.2, -0.25)}, anchor=north, legend columns=4, /tikz/every even column/.append style={column sep=0.5cm}}
    ]
    \addplot[blue, domain=0:1, samples=100] {1 - x^0.5};
    \addlegendentry{$\sigma_t=1 - t^{\alpha}$}
    \addplot[red, domain=0:1, samples=100] {(1 - x)^0.5};
    \addlegendentry{$\sigma_t=(1 - t)^{\alpha}$}
    \addplot[blue, dashed, domain=0:1, samples=100] {x^0.5};
    \addlegendentry{$\beta_t=t^{\alpha}$}
    \addplot[red, dashed, domain=0:1, samples=100] {1 - (1 - x)^0.5};
    \addlegendentry{$\beta_t=1 - (1 - t)^{\alpha}$}

    \nextgroupplot[
        xlabel={$\alpha=2$},
        width=7cm,
        height=5cm,
    ]
    \addplot[blue, domain=0:1, samples=100] {1 - x^2};
    \addplot[red, domain=0:1, samples=100] {(1 - x)^2};
    \addplot[blue, dashed, domain=0:1, samples=100] {x^2};
    \addplot[red, dashed, domain=0:1, samples=100] {1 - (1 - x)^2};
\end{groupplot}
\end{tikzpicture}
\er

We first show the following simple lemma about the density function $\phi_t(x)$.
\bl\label{Le23}
$\phi_t(x)$ solves the following first order PDE:
\begin{align}\label{PDE0}
\p_t \phi_t(x)=\left(\log\sigma_t\right)'\div ([\cD^{\mu_0}_t(x)-x]\phi_t(x)),\ \ (t,x)\in(0,T)\times\mR^d,
\end{align}
where $\cD^{\mu_0}_t(x)$ is defined by
\begin{align}\label{Def11}
\boxed{\cD^{\mu_0}_t(x):=\frac{\mE \big(\eta\rho_{\sigma_t}(x- \beta_t\eta)\big)}{\mE\rho_{ \sigma_t}(x- \beta_t\eta)}
=\frac{\mE \big(\eta\e^{-\Gamma_t(x,\eta)}\big)}{\mE\e^{-\Gamma_t(x,\eta)}}},
\end{align}
and $\eta\sim\mu_0$, and $\Gamma_t(x,y):=\|x-\beta_t y\|_2^2/(2\sigma^2_t).$
Moreover, if we let $Y_t:=\beta_t\eta+\sigma_t\xi$, where $\xi\sim N(0,\mI_d)$ is independent with $\eta$, then
\begin{align}\label{AG1}
\cD^{\mu_0}_t(x)=\mE(\eta|Y_t=x).
\end{align}
\el
\begin{proof}
First of all, recalling \eqref{DH1}, 
one sees that $\rho_{\sqrt t}(x)=(2\pi t)^{-d/2}\e^{-\|x\|^2_2/(2t)}$ is the density of a standard Brownian motion in $\mR^d$. Therefore,
$$
\p_t \rho_{\sqrt t}(x)=\tfrac12\Delta \rho_{\sqrt t}(x),\ \  t>0,\ \ x\in\mR^d.
$$
Now by the chain rule and $\beta_t=1-\sigma_t$, we have
\begin{align*}
\p_t\phi_t(x)&=\frac{(\sigma^2_t)'}{2}\int_{\mR^d}\Delta\rho_{\sigma_t}(x-\beta_ty)\mu_0(\dif y)
-\beta_t'\int_{\mR^d}y\cdot\nabla\rho_{\sigma_t}(x-\beta_ty)\mu_0(\dif y)\\
&=\sigma_t'\sigma_t\Delta \phi_t(x)+\sigma_t'\div \int_{\mR^d}y\rho_{\sigma_t}(x-\beta_ty)\mu_0(\dif y)\\
&=\div\left[\sigma_t'\sigma_t\nabla \phi_t(x)
+\sigma_t'\int_{\mR^d}y\rho_{\sigma_t}(x-\beta_ty)\mu_0(\dif y)\right].
\end{align*}
On the other hand, noting that by definition \eqref{DH1},
\begin{align*}
\nabla \phi_t(x)&=\int_{\mR^d}\frac{\beta_t y-x}{\sigma^2_t}\rho_{\sigma_t}(x-\beta_ty)\mu_0(\dif y)\\
&=\Big(\frac1{\sigma^2_t}-\frac{1}{\sigma_t}\Big)\int_{\mR^d}y\rho_{\sigma_t}(x-\beta_ty)\mu_0(\dif y)-\frac{x\phi_t(x)}{\sigma^2_t},
\end{align*}
we further have
\begin{align*}
\p_t\phi_t(x)&=\frac{\sigma'_t}{\sigma_t}\div\left[\int_{\mR^d}y\rho_{\sigma_t}(x-\beta_ty)\mu_0(\dif y)-x\phi_t(x)\right].
\end{align*}
From this we derive the desired equation \eqref{PDE0}.
For \eqref{AG1}, it is direct by definition since the density of $Y_t$ is given by $\phi_t(x)=\mE\rho_{ \sigma_t}(x- \beta_t\eta)$.
\end{proof}
\br
Compared with \eqref{PDE1},  the coefficient of PDE \eqref{PDE0} depends on the unknown distribution $\mu_0$ in an explicit way.
Moreover, for any $\epsilon_t$, we also have
\begin{align}\label{SD41}
\p_t \phi_t(x)&=\epsilon_t\Delta \phi_t(x)+\div (((\log\sigma_t)'[\cD^{\mu_0}_t(x)-x]-\epsilon_t\nabla\log \phi_t(x))\phi_t(x))\no\\
&=\epsilon_t\Delta \phi_t(x)+\div \Big(\sigma^{-2}_t\Big[\big(\sigma_t\sigma'_t-\epsilon_t\beta_t\big)\cD^{\mu_0}_t(x)+
\big(\epsilon_t-\sigma_t\sigma'_t\big)x\Big]\phi_t(x)\Big).
\end{align}
In particular, if we choose $\epsilon_t=\sigma_t\sigma'_t/\beta_t<0$ and let $p_t(x):=\phi_{1-t}(x)$ and $g_t:=-\epsilon_{1-t}$, 
$$
f_t:=\frac{\epsilon_{1-t}-\sigma_{1-t}\sigma'_{1-t}}{\sigma^2_{1-t}}=\frac{\sigma'_{1-t}}{\beta_{1-t}},
$$ 
then \eqref{SD41} becomes the forward heat equation associated with the adding noise process
$$
\p_t p_t(x)=g_t\Delta p_t(x)+\div(f_t xp_t(x)).
$$
\er

We have the following derivative estimate for nonlinear function $\cD^{\mu_0}_t(x)$ given in \eqref{Def11}.
\bl\label{Lem01}
Let $\Gamma_t(x,y):=\|x-\beta_t y\|_2^2/(2\sigma^2_t)$ and $V^{\otimes k}\in(\mR^{d})^{\otimes k}$ be the $k$-order tensor product of a vector $V\in\mR^d$.
Let $\eta\sim\mu_0$.
For any $k\in\mN$, if we write
$$
\cI_k(x):=\frac{\mE \big(\eta^{\otimes k}\e^{-\Gamma_t(x,\eta)}\big)}{\mE\e^{-\Gamma_t(x,\eta)}}\in(\mR^{d})^{\otimes k},
$$
then it holds that
$$
\nabla\cI_k(x)=\frac{\beta_t(\cI_{k+1}(x)-\cI_k(x)\otimes\cI_1(x))}{\sigma^2_t}.
$$
In particular, $\nabla\cD^{\mu_0}_t(x)\in\mR^d\otimes\mR^d$ is symmetric and positive definite, and for any $z\in\mR^d$,
\begin{align}\label{Es19}
\<z\cdot \nabla\cD^{\mu_0}_t(x), z\>=\frac{\beta_t}{\sigma^2_t}
\left[\frac{\mE \big(\<z,\eta\>^2\e^{-\Gamma_t(x,\eta)}\big)}{\mE\e^{-\Gamma_t(x,\eta)}}
-\left(\frac{\mE \big(\<z,\eta\>\e^{-\Gamma_t(x,\eta)}\big)}{\mE\e^{-\Gamma_t(x,\eta)}}\right)^2\right].
\end{align}
Moreover, if $\eta$ is bounded, then   for any $q\in[1,\infty]$ and all $t\in(0,T)$,
\begin{align}\label{Es1}
\sup_{x\in\mR^d}\|\nabla\cD^{\mu_0}_t(x)\|_q\leq 2\Big(\frac{\beta_t}{\sigma^2_t}\Big)\sup_{\omega,\omega'}\|\eta(\omega)\otimes \eta(\omega')\|_q,
\end{align}
where $\|\cdot\|_q$ stands for the $\ell^q$-norm in $(\mR^{d})^{\otimes 2}$ (see \eqref{Dis}).
\el
\begin{proof}
By the chain rule we have
\begin{align*}
\nabla\cI_k(x)
&=\frac{\mE \big(\eta^{\otimes k}\nabla\e^{-\Gamma_t(x,\eta)}\big)}{\mE\e^{-\Gamma_t(x,\eta)}}-\frac{\mE \big(\eta^{\otimes k}\e^{-\Gamma_t(x,\eta)}\big)
\otimes \mE\nabla\e^{-\Gamma_t(x,\eta)}}{(\mE\e^{-\Gamma_t(x,\eta)})^2}\\
&=\frac{\mE \big(\eta^{\otimes k}\nabla\e^{-\Gamma_t(x,\eta)}\big)}{\mE\e^{-\Gamma_t(x,\eta)}}
-\cI_k(x)\otimes\frac{\mE\nabla\e^{-\Gamma_t(x,\eta)}}{\mE\e^{-\Gamma_t(x,\eta)}}.
\end{align*}
Noting that
$$
\nabla_x\e^{-\Gamma_t(x,y)}=\left(\frac{\beta_ty-x}{\sigma^2_t}\right)\e^{-\Gamma_t(x,y)},
$$
we further have
\begin{align*}
\nabla\cI_k(x)&=\frac{\mE \big(\eta^{\otimes k}\otimes(\beta_t\eta-x)\e^{-\Gamma_t(x,\eta)}\big)}{\sigma^2_t\mE\e^{-\Gamma_t(x,\eta)}}
-\cI_k(x)\otimes\frac{\mE((\beta_t\eta-x)\e^{-\Gamma_t(x,\eta)})}{\sigma^2_t\mE\e^{-\Gamma_t(x,\eta)}}\\
&=\frac{\beta_t\mE \big(\eta^{\otimes(k+1)}\e^{-\Gamma_t(x,\eta)}\big)-\mE \big(\eta^{\otimes k}\e^{-\Gamma_t(x,\eta)}\big)\otimes x}{\sigma^2_t\mE\e^{-\Gamma_t(x,\eta)}}
-\cI_k(x)\otimes\left(\frac{\beta_t\cI_1(x)-x}{\sigma^2_t}\right)\\
&=\frac{\beta_t\cI_{k+1}(x)-\cI_k(x)\otimes x}{\sigma^2_t}-\cI_k(x)\otimes\left(\frac{\beta_t\cI_1(x)-x}{\sigma^2_t}\right)
=\frac{\beta_t(\cI_{k+1}(x)-\cI_k(x)\otimes\cI_1(x))}{\sigma^2_t}.
\end{align*}
Equality \eqref{Es19} and estimate \eqref{Es1} follow by 
$$
\nabla\cD^{\mu_0}_t(x)=\nabla\cI_1(x)=\frac{\beta_t(\cI_2(x)-\cI_1(x)\otimes\cI_1(x))}{\sigma^2_t}.
$$
The proof is complete.
\end{proof}

\subsection{Random flows of ODEs}
In this subsection we consider the random flow determined by the following ODE starting from a normal distribution:
\begin{align}\label{SDE0}
\frac{\dif Y_t}{\dif t}=(\log\sigma_t)'[Y_t-\cD^{\mu_0}_t(Y_t)]=b^{\mu_0}_t(Y_t),\ \ Y_0\sim N(0,\mI_d),
\end{align}
where
\begin{align}\label{BB0}
b^{\mu_0}_t(x):=(\log\sigma_t)'[x-\cD^{\mu_0}_t(x)].
\end{align}
We have the following simple lemma.
\bl\label{Le27}
For any $x,y\in\mR^d$, it holds that
$$
\<x-y,b^{\mu_0}_t(x)-b^{\mu_0}_t(y)\>\leq (\log\sigma_t)'\|x-y\|^2_2\left[1-\beta_t\sup_{\omega}\|\eta(\omega)\|_2^2/\sigma^2_t\right].
$$
\el
\begin{proof}
Note that 
\begin{align*}
\cD^{\mu_0}_t(x)-\cD^{\mu_0}_t(y)
=(x-y)\cdot\int^1_0\nabla\cD^{\mu_0}_t(\theta x+(1-\theta)y)\dif\theta.
\end{align*}
By \eqref{Es19}, we have
$$
\<x-y,\cD^{\mu_0}_t(x)-\cD^{\mu_0}_t(y)\>\leq \beta_t\|x-y\|^2_2\sup_{\omega}\|\eta(\omega)\|_2^2/\sigma^2_t.
$$
The desired estimate follows by $(\log\sigma_t)'<0$.
\end{proof}
The following is the main theoretical result of this section.
\bt\label{Th24}
Let $Y_0=(Y^1_0,\cdots, Y^d_0)$, where $Y^i_0\sim N(0,1)$ is a family of i.i.d. standard normal random variables. 
Let $\eta\sim\mu_0$ be independent with $Y_0$ and satisfy that  for some $q\in[1,\infty]$,
$$
\sup_{\omega}\|\eta(\omega)\|_q\leq K.
$$
Under {\bf (H)}, there is a unique solution $Y_t$ to ODE \eqref{SDE0} before reaching the terminal time $T$.
Moreover, for each $t\in[0,T)$, the law $\mu_t$ of $Y_t$ has a density given by $\phi_t(x)$, and for any $p\geq 1$,
\begin{align}\label{DQ1}
\cW_{q,p}(\mu_t,\mu_0)\leq \left[K+\left(\mE\|Y_0\|^p_q\right)^{1/p}\right]\sigma_{t},\ \ t\in[0,T),
\end{align}
and for all $t\in[0,T)$ and $z\in\mR^d$,
\begin{align}\label{UN1}
\|Y_t-\beta_tz\|_q\leq \sigma_t\|Y_0\|_q+\beta_t\sup_{\omega}\|\eta(\omega)-z\|_q.
\end{align}
In particular, for any $M>0$, 
\begin{align}\label{UN2}
\mP\left(\|Y_t\|_\infty\geq \beta_t\sup_{\omega}\|\eta(\omega)\|_\infty+\sigma_t M\right)\leq 1-(1-{\rm erfc}(M/\sqrt2))^d=:\mathfrak{p}(M,d),
\end{align}
where ${\rm erfc}(r):=\frac{2}{\sqrt\pi}\int_{r}^{\infty}\e^{-x^2}\dif x.$
\et
\begin{proof}
We divide the proof into four steps.

{\it (Step 1).} Since $\mu_0$ has compact support, by Lemma \ref{Lem01} we have
\begin{align*}
\|\nabla  b^{\mu_0}_t(x)\|_2
&\leq \left|\frac{{\sigma'_t}}{\sigma_t}\right|\left[\|\nabla\cD^{\mu_0}_{t}(x)\|_2+ \sqrt d\right]\leq C_d\left|\frac{{\sigma'_t}}{\sigma_t}\right|\left[\Big(\frac{\beta_t}{\sigma^2_t}\Big)K^2+1\right].
\end{align*}
In particular, by {\bf (H)}, for each $t_0<T$, there is a constant $C=C(d,t_0)>0$ such that
$$
\sup_{t\in[0,t_0]}\left(\|b^{\mu_0}_t(0)\|_2+\sup_{x\in\mR^d}\|\nabla  b^{\mu_0}_t(x)\|_2\right)\leq C.
$$
Hence, by the classical theory of ODE, for each starting point $x\in\mR^d$, there is a unique solution $X_t=X_t(x)$ to ODE
$$
X_t(x)=x+\int^t_0  b^{\mu_0}_s(X_s(x))\dif s,\ \ t\in[0,T).
$$
Moreover, $x\mapsto X_t(x)$ forms a smooth diffemorphism flow and the Jacobi flow $J_t(x):=\nabla X_t(x)$ satisfies
$$
J_t(x)=\mI+\int^t_0\nabla b^{\mu_0}_s(X_s(x))J_s(x)\dif s.
$$
In particular, the solution of ODE \eqref{SDE0} starting from the normal distribution $Y_0$ is given by
$$
Y_t=X_t(Y_0).
$$
Now, for any bounded measurable $f: \mR^d\to\mR$, by the change of variable, we have
$$
\mE f(Y_t)=\int_{\mR^d} f(X_t(x))\rho_1(x)\dif x
=\int_{\mR^d} f(x)\rho_1(X^{-1}_t(x))\det(\nabla X^{-1}_t(x))\dif x,
$$
which implies that $Y_t$ has a density $\psi_t(x)$ given by
$$
\psi_t(x)=\rho_1(X^{-1}_t(x))\det(\nabla X^{-1}_t(x)). 
$$

{\it (Step 2).}  By the chain rule, for any $f\in C^1_b(\mR^2)$, we also have 
$$
f(Y_t)=f(Y_0)+\int^t_0b^{\mu_0}_s(Y_s)\cdot\nabla f(Y_s)\dif s.
$$
Taking expectations and by the arbitrariness of $f$, we get
$$
\psi_t(x)=\psi_0(x)-\int^t_0\div ( b^{\mu_0}_s(x)\psi_s(x))\dif s.
$$
Let $h_t(x)=\psi_t(x)-\phi_t(x)$. Since $\psi_0=\phi_0=\rho_1$, by  \eqref{BB0} and \eqref{PDE0}, we have
$$
h_t(x)=-\int^t_0\div ( b^{\mu_0}_s(x)h_s(x))\dif s.
$$
In particular, by the integration by parts,
\begin{align*}
\p_t\int_{\mR^d}|h_t(x)|^2\dif x&=-\int_{\mR^d}\div ( b^{\mu_0}_t(x)h_t(x))h_t(x)\dif x
=\int_{\mR^d} b^{\mu_0}_t(x)h_t(x))\nabla h_t(x)\dif x\\
&=\int_{\mR^d} b^{\mu_0}_t(x)\nabla |h_t(x)|^2\dif x=-\int_{\mR^d}\div  b^{\mu_0}_t(x)|h_t(x)|^2\dif x,
\end{align*}
and
$$
\p_t\int_{\mR^d}|h_t(x)|^2\dif x\leq\|\div  b^{\mu_0}_t\|_\infty\int_{\mR^d}|h_t(x)|^2\dif x.
$$
By Gronwall's inequality, we get $\int_{\mR^d}|h_t(x)|^2\dif x=0$ for each $t\in[0,T)$. Hence,
$$
\psi_t(x)= \phi_t(x),\ \ t\in[0,T).
$$
On the other hand, by the expression \eqref{Den1},  $\phi_t(x)$ is the density of $\beta_t\eta+\sigma_t\xi,$ where $\xi\sim N(0,\mI_d)$.
Estimate \eqref{DQ1} now follows by definition \eqref{WW1} (see Lemma \ref{Le21}).

{\it (Step 3).} To show the uniform estimate \eqref{UN1}, for  $z\in\mR^d$, we introduce 
$$
Z_t(z):=(Y_t-\beta_tz)/\sigma_t.
$$
By the chain rule, we have
\begin{align*}
\frac{\dif Z_t(z)}{\dif t}=\frac{Y_t'-\beta_t'z}{\sigma_t}-\frac{(Y_t-\beta_tz)\sigma'_t}{\sigma_t^2}.
\end{align*}
Since $1-\beta_t=\sigma_t$, by  \eqref{SDE0}, we have
\begin{align*}
\frac{\dif Z_t(z)}{\dif t}&=\frac{(\log\sigma_t)'[Y_t-\cD^{\mu_0}_{t}(Y_t)]
+\sigma'_tz}{\sigma_t}-\frac{(Y_t-\beta_tz)\sigma_t'}{\sigma_t^2}
=\frac{\sigma_t'}{\sigma_t^2}[z-\cD^{\mu_0}_{t}(Y_t)].
\end{align*}
Noting that
$$
\Gamma_t(Y_t,z)=\|Y_t-\beta_tz\|_2^2/(2\sigma^2_t)=\|Z_t(z)\|_2^2/2,
$$
and recalling the definition \eqref{Def11}, we have for each $z\in\mR^d$,
\begin{align*}
\frac{\dif Z_t(z)}{\dif t}=\left(\frac{\sigma'_t}{\sigma^2_t}\right)\frac{\mE \big((z-\eta)\e^{-\Gamma_t(Y_t,\eta)}\big)}{\mE\e^{-\Gamma_t(Y_t,\eta)}}
=\left(\frac{\sigma'_t}{\sigma^2_t}\right)\frac{\mE \big((z-\eta)\e^{-\|Z_t(\eta)\|_2^2/2}\big)}{\mE \e^{-\|Z_t(\eta)\|_2^2/2}},
\end{align*}
where the expectation is taken with respect to $\eta$. 
In other words, for each $t\in[0,1)$,
\begin{align}\label{SE1}
Z_t(z)=Y_0+\int^t_0\left(\frac{\sigma_s'}{\sigma_s^2}\right)\frac{\mE \big((z-\eta) \e^{-\|Z_s(\eta)\|_2^2/2}\big)}{\mE \e^{-\|Z_s(\eta)\|_2^2/2}}\dif s.
\end{align}
Hence,
\begin{align*}
\|Z_t(z)\|_q&\leq\|Y_0\|_q+\sup_{\omega}\|\eta(\omega)-z\|_q\int^t_0\left(-\frac{\sigma_s'}{\sigma_s^2}\right)\dif s\\
&=\|Y_0\|_q+\sup_{\omega}\|\eta(\omega)-z\|_q\left(\frac1{\sigma_t}-\frac1{\sigma_0}\right),
\end{align*}
which gives \eqref{UN1}.

{\it (Step 4).} By \eqref{UN1} with $z=0$, we have
$$
\|Y_t\|_\infty\leq \sigma_t\|Y_0\|_\infty+\beta_t\sup_{\omega}\|\eta(\omega)\|_\infty.
$$
Hence, by the independence of $Y^1_0,\cdots, Y^d_0$,
\begin{align*}
&\mP\left(\|Y_t\|_\infty\geq \beta_t\sup_{\omega}\|\eta(\omega)\|_\infty+\sigma_t M\right)\leq 
\mP\left(\|Y_0\|_\infty\geq M\right)\\
&\qquad=1-\mP\left(\|Y_0\|_\infty<M\right)=1-\mP\left(|Y^1_0|<M\right)^d=1-(1-\mP\left(|Y^1_0|\geq M\right))^d.
\end{align*}
Note that
$$
\mP(|Y^1_0|\geq M)=\frac{2}{\sqrt{2\pi}}\int^\infty_M\e^{-x^2/2}\dif x=\frac{2}{\sqrt\pi}\int^\infty_{M/\sqrt2}\e^{-x^2}\dif x={\rm erfc}(M/\sqrt2).
$$
The desired estimate \eqref{UN2} now follows.
\end{proof}
\br
Estimate \eqref{DQ1} gives the weak convergence rate of $Y_t$ to $X_0$ as $t\to T$ with respect to the Wasserstein metric $\cW_{q,p}$. In particular,
$$
\cW_{2,2}(\mu_t,\mu_0)\leq \left(\sup_{\omega}\|\eta(\omega)\|_2+\sqrt{d}\right)\sigma_{t},
$$
and
$$
\cW_{\infty,2}(\mu_t,\mu_0)\leq \left(\sup_{\omega}\|\eta(\omega)\|_\infty+\sqrt{\cN(d)}\right)\sigma_{t},
$$
where
$$
\cN(\alpha):=4\int^\infty_0 r(1-(1-{\rm erfc}(r))^\alpha)\dif r.
$$
Here is the graph of function $\alpha\mapsto\sqrt{\cN(\alpha)}$. 
\begin{center}
\includegraphics[width=2.7in, height=1.5in]{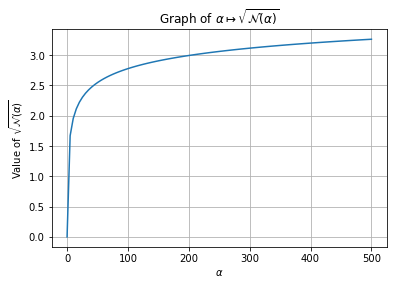}.
\includegraphics[width=2.7in, height=1.5in]{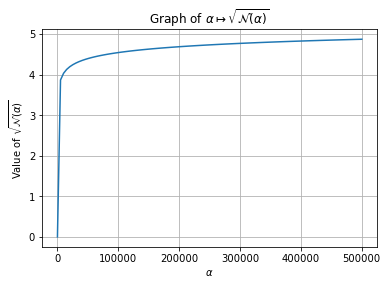}.
\end{center}
\er
\br
Some values of function $\mathfrak{p}(M,d)$ are listed below. 
\begin{center}
\begin{tabular}{|c|c|c|c|c|c|c|}
\hline
\rowcolor{gray!20}
$d$ & $\mathfrak{p}(1,d)$ & $\mathfrak{p}(2,d)$ & $\mathfrak{p}(3,d)$ & $\mathfrak{p}(4,d)$ & $\mathfrak{p}(5,d)$ & $\mathfrak{p}(6,d)$ \\
\hline
10 & 0.978010 &\cellcolor{blue!20} 0.372291 &\cellcolor{red!20} 0.026672 &\cellcolor{green!20}0.000633
 & 0.000006 & 0.000000 \\
\hline
100 & 1.000000 & 0.990503 &\cellcolor{blue!20} 0.236884 & \cellcolor{red!20}0.006314 & \cellcolor{green!20}0.000057 & 0.000000 \\
\hline
1000 & 1.000000 & 1.000000 &\cellcolor{blue!20} 0.933026 &\cellcolor{red!20} 0.061380 &\cellcolor{green!20} 0.000573 & 0.000002 \\
\hline
10000 & 1.000000 & 1.000000 & 1.000000 &\cellcolor{blue!20} 0.469240 &\cellcolor{red!20}  0.005717 & \cellcolor{green!20}  0.000020 \\
\hline
100000 & 1.000000 & 1.000000 & 1.000000 &\cellcolor{blue!20} 0.998226 &\cellcolor{red!20} 0.055718 &\cellcolor{green!20} 0.000197 \\
\hline
\end{tabular}
\end{center}
From the above table, one sees that for $d=10000$,
$$
\mP\left(\|Y_t\|_\infty\geq \beta_t\sup_{\omega}\|\eta(\omega)\|_\infty+5\sigma_t\right)\leq 0.005717.
$$
\er

Next starting from \eqref{SE1}, we derive an useful equation for $\mE\e^{-\|Z_t(\eta)\|_2^2/2}$.
\bt
Let $g_t:=\mE\e^{-\Gamma_t(Y_t,\eta)}$, where $\Gamma_t(x,y):=\|x-\beta_t y\|_2^2/(2\sigma^2_t)$ and the expectation is only taken with respect to $\eta$, not for $Y_t$. 
Then $g_0=\e^{-\|Y_0\|_2^2/2}$ and it holds that
\begin{align}\label{Eq1}
g_t'=\left(\frac{\sigma_t'\beta_t}{\sigma^3_t}\right)\left[\Lambda_2(t,Y_t)-\frac{\|\Lambda_1(t,Y_t)\|^2_2}{g_t}\right],
\end{align}
where
\begin{align}\label{No1}
\Lambda_1(t,x):=\mE\left[\eta\e^{-\Gamma_t(x,\eta)}\right],\ \ \Lambda_2(t,x):=\mE\left[\|\eta\|^2_2\e^{-\Gamma_t(x,\eta)}\right].
\end{align}
In particular, suppose that $\sup_\omega\|\eta(\omega)\|_2\leq K$, then
\begin{align}\label{Bo1}
\e^{-[\|Y_0\|_2^2+(K\beta_t/\sigma_t)^2]/2}\leq g_t\leq \e^{-\|Y_0\|_2^2/2},\ \ t\in[0,T).
\end{align}
\et
\begin{proof}
Recall $Z_t(z):=(Y_t-\beta_tz)/\sigma_t.$
By \eqref{SE1} and the chain rule, we  have for each $z\in\mR^d$,
$$
\|Z_t(z)\|_2^2/2=\|Y_0\|_2^2/2+\int^t_0\left(\frac{\sigma_s'}{\sigma_s^2}\right)
\frac{\mE\big(\<Z_s(z),z-\eta\> \e^{-\|Z_s(\eta)\|_2^2/2}\big)}{\mE \e^{-\|Z_s(\eta)\|_2^2/2}}\dif s,
$$
and
$$
\e^{-\|Z_t(z)\|_2^2/2}=\e^{-\|Y_0\|_2^2/2}-\int^t_0\left(\frac{\sigma_s'}{\sigma_s^2}\right)
\frac{\mE\big(\<Z_s(z),z-\eta\> \e^{-(\|Z_s(z)\|_2^2+\|Z_s(\eta)\|_2^2)/2}\big)}{\mE \e^{-\|Z_s(\eta)\|_2^2/2}}\dif s.
$$
Let $\tilde\eta$ be an independent copy of $\eta$, which is also independent with $Y_0$. Then replacing the above $z$ by $\tilde\eta$ and taking expectations with respect to $\tilde\eta$, we obtain
\begin{align}\label{HF1}
g_t=\mE\e^{-\|Z_t(\tilde\eta)\|_2^2/2}=\e^{-\|Y_0\|_2^2/2}-\int^t_0\left(\frac{\sigma_s'}{\sigma_s^2}\right)
\frac{\mE\big(\<Z_s(\tilde\eta),\tilde\eta-\eta\> \e^{-(\|Z_s(\tilde\eta)\|_2^2+\|Z_s(\eta)\|_2^2)/2}\big)}{g_s}\dif s.
\end{align}
Recalling that $Z_t(\tilde\eta):=(Y_t-\beta_t\tilde\eta)/\sigma_t$ and $\eta,\tilde\eta$ are independent, we have
\begin{align*}
&\mE\big(\<Z_s(\tilde\eta),\tilde\eta-\eta\> \e^{-(\|Z_s(\tilde\eta)\|_2^2+\|Z_s(\eta)\|_2^2)/2}\big)\\
&=\frac1{\sigma_s}\mE\left[\<Y_s-\beta_s\tilde\eta, \tilde\eta-\eta\>\e^{-(\|Z_s(\tilde\eta)\|_2^2+\|Z_s(\eta)\|_2^2)/2} \right]\\
&=\frac{\beta_s}{\sigma_s}\mE\left[\<\tilde\eta, \eta-\tilde\eta\>\e^{-(\|Z_s(\tilde\eta)\|_2^2+\|Z_s(\eta)\|_2^2)/2} \right]\\
&=\frac{\beta_s}{\sigma_s}\left[\|\Lambda_1(t,Y_s)\|^2_2-\Lambda_2(t,Y_s)g_s\right],
\end{align*}
where we have used that (the expectations are taken only with respect to $\tilde\eta$ and $\eta$)
$$
\mE\left[\<Y_s, \tilde\eta-\eta\> \e^{-(\|Z_s(\tilde\eta)\|_2^2+\|Z_s(\eta)\|_2^2)/2}\right]=0.
$$
Substituting it into \eqref{HF1}, we obtain  equation \eqref{Eq1}.

To prove \eqref{Bo1}, note that by H\"older's inequality,
$$
\|\Lambda_1(t,Y_t)\|^2_2\leq g_t\Lambda_2(t,Y_t).
$$
Since $\sigma'_t<0$ and $\sup_\omega\|\eta(\omega)\|_2\leq K$, we have
$$
\left(\frac{\sigma_t'\beta_t}{\sigma^3_t}\right)K^2\, g_t\leq g'_t\leq 0.
$$
Solving this differential inequality, we derive that
$$
\e^{K^2\int^t_0 \sigma_s'\beta_s/\sigma^3_s\dif s}g_0\leq g_t\leq g_0,\ \ t\in(0,T).
$$
Since $\beta_t=1-\sigma_t$, by the change of variable $\sigma_s=r$, it is easy to see that
\begin{align}\label{SUP1}
\int^t_0 \sigma_s'\beta_s/\sigma^3_s\dif s=\int^{\sigma_0}_{\sigma_t} (1-r)/r^3\dif r=-\frac{\beta^2_t}{2\sigma_t^2}.
\end{align}
The proof is complete.
\end{proof}
\br
By using the  definitions \eqref{No1}, we also have
$$
Y'_t=\left(\frac{\sigma_t' }{\sigma_t}\right)\left[Y_t-\frac{\Lambda_1(t,Y_t)}{g_t}\right].
$$
In other words, $(Y_t, g_t)$ solves the following ODE system:
$$
\left\{
\begin{aligned}
Y_t&=Y_0+\int^t_0\left(\frac{\sigma_s' }{\sigma_s}\right)\left[Y_s-\frac{\Lambda_1(s,Y_s)}{g_s}\right]\dif s,\\
g_t&=\e^{-\|Y_0\|_2^2/2}+\int^t_0
\left(\frac{\sigma_s'\beta_s}{\sigma^3_s}\right)\left[\Lambda_2(s,Y_s)-\frac{\|\Lambda_1(s,Y_s)\|^2_2}{g_s}\right]\dif s.
\end{aligned}
\right.
$$
\er

\section{Application to diffusion generative models}

In this section we aim to introduce an algorithm for generating a sample from the dataset 
$$
\cX_0=\{\eta^1_0,\cdots, \eta^N_0\},
$$ 
where each data $\eta^j_0$ in $\cX_0$ obeys the same distribution $\mu_0$. We suppose that
\begin{align}\label{DS1}
\mbox{$\mu_0$ has support in $\big\{x\in\mR^d: \|x\|_2\leq K\big\}$}.
\end{align}
Let $\mu^N_0$ be the empirical measure of $\cX_0$, i.e.,
\begin{align}\label{Def01}
\mu^N_0(\dif y):=\frac1N\sum_{j=1}^N\delta_{\eta^j_0}(\dif y).
\end{align}
Consider the particle approximation coefficients:
\begin{align}\label{Ex19}
\cD^{\mu_0^N}_t(x)
&=\frac{\int_{\mR^d}y\rho_{\sigma_t}(x-\beta_t y)\mu^N_0(\dif y)}{\int_{\mR^d}\rho_{\sigma_t}(x-\beta_t y)\mu^N_0(\dif y)}.
\end{align}
\bl\label{Le51}
Let $\phi_t(x)$ be defined by \eqref{Den1} and $\cD^{\mu_0}_t(x)$ be defined by \eqref{Def11}. 
We have
$$
\|\cD^{\mu_0}_t(x)-\cD^{\mu^N_0}_t(x)\|_2\leq \sI^N_t(x)/\phi_t(x),
$$
where 
\begin{align}\label{AG2}
\sI^N_t(x):=\left\|\int_{\mR^d}y\rho_{\sigma_t}(x-\beta_t y)(\mu_0-\mu^N_0)(\dif y)\right\|_2+K\left|\int_{\mR^d}\rho_{\sigma_t}(x-\beta_t y)(\mu_0-\mu^N_0)(\dif y)\right|.
\end{align}
Moreover, it holds that
\begin{align}\label{AG20}
\|\nabla\cD^{\mu_0^N}_t(x)\|_2\leq 2\beta_tK^2/\sigma^2_t.
\end{align}
\el
\begin{proof}
By definition we have
\begin{align*}
\|\cD^{\mu_0}_t(x)-\cD^{\mu^N_0}_t(x)\|_2&=\left\|\frac{\int_{\mR^d}y\rho_{\sigma_t}(x-\beta_t y)\mu_0(\dif y)}{\int_{\mR^d}\rho_{\sigma_t}(x-\beta_t y)\mu_0(\dif y)}
-\frac{\int_{\mR^d}y\rho_{\sigma_t}(x-\beta_t y)\mu^N_0(\dif y)}{\int_{\mR^d}\rho_{\sigma_t}(x-\beta_t y)\mu^N_0(\dif y)}\right\|_2\\
&\leq\frac{\left\|\int_{\mR^d}\rho_{\sigma_t}(x-\beta_t y)\mu^N_0(\dif y)\int_{\mR^d}y\rho_{\sigma_t}(x-\beta_t y)(\mu_0-\mu^N_0)(\dif y)\right\|_2}{\int_{\mR^d}\rho_{\sigma_t}(x-\beta_t y)\mu_0(\dif y)\int_{\mR^d}\rho_{\sigma_t}(x-\beta_t y)\mu^N_0(\dif y)}\\
&+\frac{\left\|\int_{\mR^d}y\rho_{\sigma_t}(x-\beta_t y)\mu^N_0(\dif y)\int_{\mR^d}\rho_{\sigma_t}(x-\beta_t y)(\mu_0-\mu^N_0)(\dif y)\right\|_2}{\int_{\mR^d}\rho_{\sigma_t}(x-\beta_t y)\mu_0(\dif y)\int_{\mR^d}\rho_{\sigma_t}(x-\beta_t y)\mu^N_0(\dif y)}\\
&\leq\frac{\left\|\int_{\mR^d}y\rho_{\sigma_t}(x-\beta_t y)(\mu_0-\mu^N_0)(\dif y)\right\|_2}{\int_{\mR^d}\rho_{\sigma_t}(x-\beta_t y)\mu_0(\dif y)}+\frac{K\left|\int_{\mR^d}\rho_{\sigma_t}(x-\beta_t y)(\mu_0-\mu^N_0)(\dif y)\right|}{\int_{\mR^d}\rho_{\sigma_t}(x-\beta_t y)\mu_0(\dif y)},
\end{align*}
which gives the first estimate. As for \eqref{AG20}, it follows by \eqref{Es1}.
\end{proof}
Let $Y_0\sim N(0, \mI_d)$ be independent of $\cX_0$.
We consider the following approximation ODE:
\begin{align}\label{ODE0}
\frac{\dif Y^N_t}{\dif t}=\left(\frac{{\sigma'_t}}{\sigma_t}\right) [Y^N_t-\cD^{\mu^N_0}_t(Y^N_t)]=b^N_t(Y^N_t),\ \ Y^N_0=Y_0.
\end{align}
We have the following convergence estimate.
\bt\label{Th51}
Under {\bf (H)} and \eqref{DS1}, it holds that for all $t\in[0,T)$,
$$
\sup_{N\in\mN}\left(\sqrt N\mE\|Y_t-Y^N_t\|_2\right)\leq 2K\beta_t\e^{K^2\beta_t^2/(2\sigma^2_t)},
$$
where the expectation is taken with respect to $\cX_0$ and $Y_0$.
\et
\begin{proof}
Note that
$$
Y_t-Y^N_t=\int^t_0 ( b_s(Y_s)- b^N_s(Y^N_s))\dif s
=\int^t_0 ( b_s(Y_s)- b^N_s(Y_s))\dif s+\int^t_0 ( b^N_s(Y_s)- b^N_s(Y^N_s))\dif s.
$$
By the chain rule and Lemmas \ref{Le27} and \ref{Le51}, we have
\begin{align*}
\frac{\dif \|Y_t-Y^N_t\|_2}{\dif t}
&=\frac{\<Y_t-Y^N_t, b_t(Y_t)- b^N_t(Y_t)\>+\<Y_t-Y^N_t, b^N_t(Y_t)- b^N_t(Y^N_t)\>}{\|Y_t-Y^N_t\|_2}\\
&\leq\frac{|\sigma_t'|}{\sigma_t}\left(\|\cD^{\mu_0}_t(Y_t)-\cD^{\mu^N_0}_t(Y_t)\|_2+ \|Y_t-Y^N_t\|_2\left[K^2\beta_t/\sigma^2_t-1\right]\right)\\
&\leq\frac{|\sigma_t'|}{\sigma_t}\left(\frac{\sI^N_t(Y_t)}{\phi_t(Y_t)}+ \|Y_t-Y^N_t\|_2\left[K^2\beta_t/\sigma^2_t-1\right]\right).
\end{align*}
Note that by \eqref{SUP1},
$$
\int^t_0\frac{|\sigma_s'|}{\sigma_s}\left[K^2\beta_s/\sigma^2_s-1\right]\dif s=K^2\beta^2_t/(2\sigma_t^2)+\log\sigma_t.
$$
By Gronwall's inequality, we get
\begin{align*}
\|Y_t-Y^N_t\|_2&\leq\int^t_0\e^{\int^t_s\frac{|\sigma_r'|}{\sigma_r}\left[K^2\beta_r/\sigma^2_r-1\right]\dif r}
\frac{|\sigma_s'|}{\sigma_s}\frac{\sI^N_s(Y_s)}{\phi_s(Y_s)}\dif s\\
&=\e^{K^2\beta^2_t/(2\sigma_t^2)}\sigma_t
\int^t_0\e^{-K^2\beta^2_s/(2\sigma_s^2)}\frac{|\sigma_s'|}{\sigma^2_s}\frac{\sI^N_s(Y_s)}{\phi_s(Y_s)}\dif s,
\end{align*}
which implies by taking expectations that
\begin{align}\label{DQ1}
\mE\|Y_t-Y^N_t\|_2&\leq\e^{K^2\beta^2_t/(2\sigma_t^2)}\sigma_t
\int^t_0\e^{-K^2\beta^2_s/(2\sigma_s^2)}\frac{|\sigma_s'|}{\sigma^2_s}\mE\left[\frac{\sI^N_s(Y_s)}{\phi_s(Y_s)}\right]\dif s.
\end{align}
Noting that $Y_s$ is independent with $\cX_0$ and has density $\phi_s(x)$, we   have
\begin{align*}
\mE\left[\frac{\sI^N_s(Y_s)}{\phi_s(Y_s)}\right]=\int_{\mR^d}\mE^{\cX_0}[\sI^N_s(x)]\dif x,
\end{align*}
where $\mE^{\cX_0}$ stands for the expectation with respect to $\cX_0$.
Since for a function $\varphi:\mR^d\to\mR^m$,
$$
\mE^{\cX_0}\left\|\int_{\mR^d}\varphi(y)(\mu_0-\mu^N_0)(\dif y)\right\|_2^2
=\frac1N\left(\mE\|\varphi(\eta)\|^2_2-\left\|\mE\varphi(\eta)\right\|^2_2\right),
$$
where $\eta\sim\mu_0$ is independent with $\cX_0$ and $Y_0$, by \eqref{DS1},
we have
\begin{align*}
\mE^{\cX_0}\left\|\int_{\mR^d}y\rho_{\sigma_s}(x-\beta_sy)(\mu_0-\mu^N_0)(\dif y)\right\|_2^2
&\leq\frac{\mE(\|\eta\|^2_2\rho^2_{\sigma_s}(x-\beta_s\eta))}N
\leq \frac{K^2\mE\rho^2_{\sigma_s}(x-\beta_s\eta)}N,
\end{align*}
and also,
\begin{align*}
\mE^{\cX_0}\left|\int_{\mR^d}\rho_{\sigma_s}(x-\beta_sy)(\mu_0-\mu^N_0)(\dif y)\right|^2
&\leq \frac{\mE\rho^2_{\sigma_s}(x-\beta_s\eta)}N.
\end{align*}
Hence, by H\"older's inequality and \eqref{AG2}, 
$$
\mE^{\cX_0}[\sI^N_s(x)]\leq \Big(\mE^{\cX_0}[\sI^N_s(x)]^2\Big)^{1/2}\leq\frac{2K}{\sqrt N}\left(\mE\rho^2_{\sigma_s}(x-\beta_s\eta)\right)^{1/2},
$$
and by Jensen's inequality and recalling the expression \eqref{DH1},
\begin{align*}
\mE\left[\frac{\sI^N_s(Y_s)}{\phi_s(Y_s)}\right]&
=\int_{\mR^d}\mE^{\cX_0}[\sI^N_s(x)]\dif x\leq\frac{2K}{\sqrt N}\int_{\mR^d}\left(\mE\rho^2_{\sigma_s}(x-\beta_s\eta)\right)^{1/2}\dif x\\
&\leq\frac{2K}{\sqrt N}\left(\int_{\mR^d}\mE\rho^2_{\sigma_s}(x-\beta_s\eta) \rho^{-1}_{\sigma_s}(x)\dif x\right)^{1/2}\\
&=\frac{2K}{\sqrt N}\left(\mE\int_{\mR^d}\frac{\e^{(\|x\|^2_2-2\|x-\beta_s \eta\|^2_2)/(2\sigma^2_s)}}{(2\pi\sigma^2_s)^{d/2}}\dif x\right)^{1/2}\\
&=\frac{2K}{\sqrt N}\left(\mE\int_{\mR^d}\frac{\e^{(2\|\beta_s\eta\|^2_2-\|x-2\beta_s \eta\|^2_2)/(2\sigma^2_s)}}{(2\pi\sigma^2_s)^{d/2}}\dif x\right)^{1/2}\\
&=\frac{2K}{\sqrt N}\Big(\mE\e^{\|\beta_s\eta\|^2_2/\sigma^2_s}\Big)^{1/2}\leq\frac{2K}{\sqrt N}\e^{K^2\beta^2_s/(2\sigma^2_s)}.
\end{align*}
Substituting this into \eqref{DQ1} and by $\beta_s=1-\sigma_s$, we obtain
\begin{align*}
\mE\|Y_t-Y^N_t\|_2&\leq\frac{2K\e^{K^2\beta^2_t/(2\sigma_t^2)}\sigma_t}{\sqrt N}\int^t_0\frac{|\sigma_s'|}{\sigma_s^2}\dif s
=\frac{2K\e^{K^2\beta^2_t/(2\sigma_t^2)}\beta_t}{\sqrt N}.
\end{align*}
The proof is complete.
\end{proof}
\br
It must be noted that our convergence estimate does not depend on the dimension.
\er

Recalling $\Gamma_t(x,y):=-\|x-\beta y\|^2_2/(2\sigma^2_t)$,
we can write \eqref{Ex19} as
$$
\cD^{\mu_0^N}_t(x)
=\frac{\int_{\mR^d}y\e^{-\Gamma_t(x,y)}\mu^N_0(\dif y)}{\int_{\mR^d}\e^{-\Gamma_t(x,y)}\mu^N_0(\dif y)}
=\frac{\sum_{j=1}^N \big(\eta^j_0\e^{-\Gamma_t(x,\eta_0^j)}\big)}{\sum_{j=1}^N\e^{-\Gamma_t(x,\eta_0^j)}}.
$$
Using this expression
and the classical Euler scheme for ODE, we have the following algorithm for generating a sample from the dataset $\cX_0$.

\noindent\rule{\linewidth}{1 pt}

{\bf Algorithm 1.} Generate sample from data without learning

\noindent\rule{\linewidth}{0.8 pt}

\begin{enumerate}[1.]
\item {\bf Data}: Sample data $\{\eta^j_0\}_{j=1}^N$, time iteration number $M$.
\item {\bf Initial values}: Generate normal random variables $\xi\sim N(0, \mI_d)$.
\item {\bf Normalizing initial values:} $Y_0=\xi \sqrt{d}/\|\xi\|_2$ so that $\|Y_0\|^2_2=d$. 
\item {\bf For} $k=0$ {\bf to} $M-1$
\begin{itemize}
\item $\sigma_{k}=1-k/M$; $\beta_k=k/M$;
\item $Y_{k+1}=Y_k+\left[\frac{\sum_{j=1}^N \big(\eta^j_0\e^{d/2-\|Y_k-\beta_k\eta^j_0\|_2^2/(2\sigma^2_k)}\big)}{\sum_{j=1}^N\e^{d/2-\|Y_k-\beta_k\eta^j_0\|_2^2/(2\sigma^2_k)}}-Y_k\right]/(M-k)$.
\end{itemize}
\item {\bf Output}: $Y_M$ as generated sample.
\end{enumerate}
\noindent\rule{\linewidth}{1 pt}

\br
Here, when $d$ is large, to avoid the overflow problem of the exponential function, we normalize the initial value 
$Y_0$ so that $\|Y_0\|^2_2=d$. This normalization is essential to ensure numerical stability and to maintain meaningful comparisons across different dimensions. For more details, see equation \eqref{Bo1}.
\er

Here are the results of our numerical experiment. For a given dimension $d$ and the size $N$ of dataset, 
we generate  $N$ i.i.d. uniformly distributed random variables on $[0,1]^d$ 
as our experimental dataset $\cX_0=\{\eta^j_0\}_{j=1}^N$. For different time iteration number $M$, 
we then compare the difference between the output $Y_M$ after 
$M$ iterations and the dataset $\cX_0$. The comparison metric used is the 
$L^1$-norm, defined as:
$$
\|Y_M-\cX_0\|_1=\min_{j=1,\cdots,N}\|Y_M-\eta^j_0\|_1.
$$

Below, we fix $N=10,000$ and consider different dimensions $d$ and iteration counts $M$.
The details of the numerical experiment results are summarized in the following table:

\begin{center}
\begin{tabular}{|c|c|c|c|c|c|c|c|}
\hline
 & $M=2$ & $M=3$ & $M=10$ & $M=100$ & $M=500$ & $M=1000$ & $M=2000$ \\
\hline
$d=2$ & $0.0056$ & $0.0043$ & 0.0050 &0.0054 & 0.0005 & 1.7368e-13 & 4.1633e-17 \\
\hline
$d=10$ & $0.6251$ & $0.8854$ & 0.6336 &2.151e-16 & 2.255e-16 & 2.77e-16 & -- \\
\hline
$d=100$ & $18.40325$ & $18.4920$ & 2.48e-15 &4.211e-15 & 1.576e-15 & -- & -- \\
\hline
$d=1000$ & $86.3443$ & 8.228e-14 & 4.149e-14 &2.607e-14 & -- & -- & -- \\
\hline
$d=10000$ & $755.12644$ & 7.510e-13 & 3.189e-13 &-- & -- & -- & -- \\
\hline
$d=50000$ & $24703.796$ & 4.221e-12 & -- &-- & -- & -- & -- \\
\hline
\end{tabular}
\end{center}

From the above table, it is interesting to notice that for larger dimensions 
$d$, fewer iteration times are required to obtain convergence. 
At this stage, I do not have a complete explanation for this phenomenon. 
However, a possible explanation is that in higher dimension space, the $10,000$-points are sparsely and  uniformly distributed in the space.
The Euclidean distance between each two points is comparable to the dimension $d$ (see the column of $M=2$),
which possible leads to the faster convergence. 
Moreover, I would like to highlight that the factor $d/2$ in the exponential plays an important role.
For image generation tasks, the dimension
$d$ is usually large. Hence, one can utilize the above algorithm to generate training data and then employ a neural network to learn the function 
$b^{\mu_0}_t(x)$ for generating new samples. 
We will further investigate this in future work.

 \section{Application to sampling from a known distribution}\label{Sec3}

In this section we use ODE \eqref{SDE0} to introduce a new algorithm for sampling data from a known distribution 
$\mu_0(\dif x)\varpropto f(x)\dif x$, where
\begin{align}\label{DX1}
\mbox{$f:\mR^d\to[0,\infty)$ is bounded by $\Lambda$ and has support in $\{x: \|x\|_2\leq K\}$.}
\end{align}
We emphasize that  unlike the classical Langevin dynamic method, 
we do not make any regularity assumption about $f$. In this case, we need to calculate the function 
$\cD^f_t(x):=\cD^{\mu_0}_t(x)$, that is,
\begin{align}\label{DX91}
\cD^f_t(x)=\frac{\int_{\mR^d}y\rho_{\sigma_t}(x-\beta_ty)f(y)\dif y}{\int_{\mR^d}\rho_{\sigma_t}(x-\beta_ty)f(y)\dif y}
=\frac{\int_{\mR^d}y\e^{-\Gamma_t(x,y)}f(y)\dif y}{\int_{\mR^d}\e^{-\Gamma_t(x,y)}f(y)\dif y},
\end{align}
where $\rho_\sigma$ is defined by \eqref{DH1} and $\Gamma_t(x,y)=\|x-\beta_t y\|_2^2/(2\sigma^2_t)$.
When the dimension $d$ is small, the integral $\int_{\mR^d}\e^{-\Gamma_t(x,y)}f(y)\dif y$ can be calculated by the classical numerical method. 
When $d$ is large, the traditional numerical method becomes inefficient. However, we can use the Monte-Carlo method to
approximate $\cD^f_t(x)$. 

Fix $N\in\mN$. Let $\Xi:=\{\xi_j,j=1,\cdots, N\}$ be 
a sequence of i.i.d random variables with common uniformly distribution in the ball $\{x: \|x\|_2\leq K\}$. 
Consider the approximation coefficients:
$$
\cD^f_{N,t}(x)
=\frac{\sum_{j=1}^N \xi_j\rho_{\sigma_t}(x-\beta_t\xi_j) f(\xi_j)}{\sum_{j=1}^N\rho_{\sigma_t}(x-\beta_t\xi_j)f(\xi_j)}
=\frac{\sum_{j=1}^N \xi_j\e^{-\Gamma_t(x,\xi_j)} f(\xi_j)}{\sum_{j=1}^N\e^{-\Gamma_t(x,\xi_j)} f(\xi_j)}.
$$
The following lemma is completely the same as Lemma \ref{Le51}.
\bl\label{Le61}
Let $\eta$ be the uniformly distributed random variable in $\{x: \|x\|_2\leq K\}$ and
$$
\phi_t(x):=\int_{\mR^d}\rho_{\sigma_t}(x-\beta_ty)f(y)\dif y=\mE(\rho_{\sigma_t}(x-\beta_t\eta)f(\eta)).
$$ 
We have
$$
\|\nabla\cD^f_{N,t}(x)\|_2\leq 2\beta_tK^2/\sigma^2_t,\ \ \|\cD^f_t(x)-\cD^f_{N,t}(x)\|_2\leq \sI^N_t(x)/\phi_t(x),
$$
where 
\begin{align}\label{AG2}
\begin{split}
\sI^N_t(x)&:=\left\|\frac1N\sum_{j=1}^N\xi_j\rho_{\sigma_t}(x-\beta_t \xi_j)f(\xi_j)-\mE(\eta\rho_{\sigma_t}(x-\beta_ty)f(\eta))\right\|_2\\
&\quad+K\left|\frac1N\sum_{j=1}^N\rho_{\sigma_t}(x-\beta_t \xi_j)f(\xi_j)-\mE(\rho_{\sigma_t}(x-\beta_t\eta)f(\eta))\right|.
\end{split}
\end{align}
\el

Now we can consider the following approximation ODE:
\begin{align}\label{SDEA}
\frac{\dif Y^N_t}{\dif t}=(\log\sigma_t)'\left[Y^N_t-\cD^f_{N,t}(Y^N_t)\right],\ \ Y^N_0\sim N(0,\mI_d).
\end{align}
By Lemma \ref{Le61} and as in the proof of Theorem \ref{Th51}, we also have the following dimension-free convergence result.
\bt\label{Th61}
Under {\bf (H)} and \eqref{DX1}, it holds that for all $t\in[0,1)$,
\begin{align}\label{DB1}
\sup_{N\in\mN}\left(\sqrt N\mE\|Y_t-Y^N_t\|_2\right)\leq 2K\Lambda\beta_t\e^{K^2\beta_t^2/(2\sigma^2_t)},
\end{align}
where the expectation is taken with respect to $\Xi$ and $Y_0$.
\et
\br
For larger $K$, by \eqref{DB1}, the convergence becomes worse. However, we can use the scaling technique to overcome this difficulty. Indeed, if we fix a scaling parameter $\epsilon \in (0, 1)$ and let 
$$ 
f_{K,\epsilon}(x) := f\left(Kx/\epsilon\right), 
$$ 
then $ f_{K,\epsilon} $ has support in the ball $\{x: \|x\|_2 \leq \epsilon\}$. Thus, one can replace $ K $ 
in \eqref{DB1} with $\epsilon$.
\er

It remains to discretize ODE \eqref{SDEA} by Euler's scheme. Since $x\mapsto  \cD^f_{N,t}(x)$ is smooth (see Lemma \ref{Lem01}) and there are numerous references about the convergence analysis of Euler's scheme, we will not discuss this topic.
Here is the concrete algorithm  for sampling by taking $\sigma_t=1-t$.  

\noindent\rule{\linewidth}{1 pt}

{\bf Algorithm 2.} Generate samples from known distributions.

\noindent\rule{\linewidth}{0.8 pt}

\begin{enumerate}[1.]
\item {\bf Data}: $\mu_0(\dif x)=f(x)\dif x$ with $f$ being a density function in $B_K:=\{x: \|x\|_2\leq K\}$. 
\item {\bf Parameters}: $M$ iteration number, $n$ the number of Monte-Carlo sampling, $\epsilon$ scaling parameter.
\item {\bf Initial values}: Generate a normal random variable $Y_0\sim N(0, \mI_d)$.
\item {\bf Iteration:} for $k=0$ to $M-1$
\begin{itemize}
\item Generate i.i.d. uniformly distributed random variables $\xi_1,\cdots,\xi_n$ in $B_\epsilon$.
\item $\sigma_{k}=1-k/M$; $\beta_k=k/M$.
\item $Y_{k+1}=Y_k+\left[\frac{\sum_{j=1}^N \big(\xi_j\e^{d/2-\|Y_k-\beta_k\xi_j\|_2^2/(2\sigma^2_k)}f(K\xi_j/\epsilon)\big)}{\sum_{j=1}^N\big(\e^{d/2-\|Y_k-\beta_k\xi_j\|_2^2/(2\sigma^2_k)}f(K\xi_j/\epsilon)\big)}-Y_k\right]/(M-k)$.
\end{itemize}
\item {\bf Output}: a sample $Y_M/\eps\sim\mu_0$.
\end{enumerate}
\noindent\rule{\linewidth}{1 pt}

Next, we present several numerical experimental results. In the first example, we sample from four one dimensional probability density functions, 
where the first density function is discontinuous and the others are continuous including the semicircle law. 
In the second example, we sample from 
four two-dimensional multimodal density functions. 
The figures show that the sampling results effectively simulate the distributions.
In the third example, we consider the high dimensional anisotropic funnel distribution.

{\bf Example 1.}  We consider the following one-dimensional probability density functions:
\begin{align*}
f_1(x) &= \frac{1.2 \exp(-2x^2) \cdot \mathbf{1}_{x \leq 0.5} + 2 \exp\left(-\frac{(x-1)^2}{8}\right) \cdot \mathbf{1}_{x > 0.5}}{2.8996 \sqrt{2\pi}},\\
f_2(x)&=\{200*[(x-0.7)\1_{x\in(0.7,0.8)}+(0.9-x)\1_{x\in(0.8,0.9)}]\\
&\quad+50*[(x-0.4)\1_{x\in(0.4,0.5)}+(0.6-x)\1_{x\in(0.5,0.6)}]\\
&\quad+50*[(x-0.1)\1_{x\in(0.1,0.2)}+(0.3-x)\1_{x\in(0.2,0.3)}]\}/3,
\end{align*}
and
\begin{align*}
f_3(x)=\tfrac2\pi\sqrt{1-y^2}\1_{|y|\leq 1},\ \ f_4(x)&=(1+(\sin (2\pi x)+\sin (4\pi x))/2)\1_{x\in[0,1]}.
\end{align*}
For the mentioned density functions, we generated $20,000$ sample points using the algorithm described above and plotted their histograms in Figure \ref{fig:1D} below. We used $M=50$ and 
$n=20,000$ for calculating the integral via the Monte Carlo method. The red curve represents the true distribution, while the sky-blue bars depict the histogram.

\begin{figure}[ht]
    \centering
    \begin{subfigure}[b]{0.4\textwidth}
        \includegraphics[width=2.4in, height=1.4in]{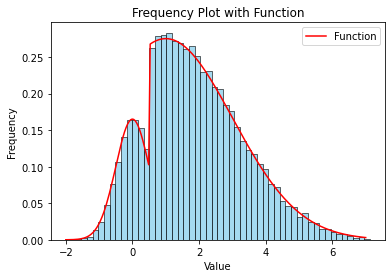}
        \caption{$f_1(x)$}
        \label{fig:f1}
    \end{subfigure}
    \hfill
    \begin{subfigure}[b]{0.4\textwidth}
        \includegraphics[width=2.4in, height=1.4in]{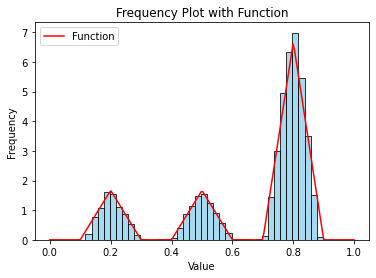}
        \caption{$f_2(x)$}
        \label{fig:f2}
    \end{subfigure}
    \begin{subfigure}[b]{0.4\textwidth}
        \includegraphics[width=2.4in, height=1.4in]{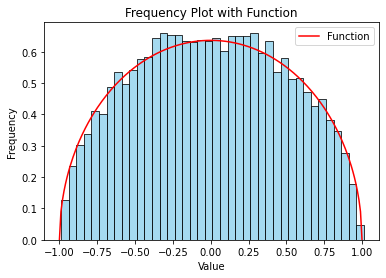}
        \caption{$f_3(x)$}
        \label{fig:f3}
    \end{subfigure}
    \hfill
    \begin{subfigure}[b]{0.4\textwidth}
        \includegraphics[width=2.4in, height=1.4in]{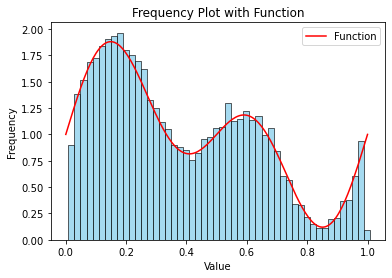}
        \caption{$f_4(x)$}
        \label{fig:f4}
    \end{subfigure}
 \caption{Sampling of 1D-probability density functions}
    \label{fig:1D}
\end{figure}

{\bf Example 2.}  Consider the following four $2$D density functions
$$
f_1(x)=c_1[(x_1^2 + x_2^2) / 4000-\cos(x_1)\cos(x_2/\sqrt 2)+1],\ \ x=(x_1,x_2)\in[0,1]^2,
$$
$$
f_2(x)=c_2\left[\e^{-\frac{\|x-(0.8,0.8)\|_2^2}{0.01}}+\e^{-\frac{\|x-(0.2,0.2)\|_2^2}{0.01}}+\e^{-\frac{\|x-(0.8,0.2)\|_2^2}{0.01}}
+\e^{-\frac{\|x-(0.2,0.8)\|_2^2}{0.01}}\right],\ \ x\in[0,1]^2,
$$ 
and
$$
f_3(x)=c_3\e^{-((x_1x_2)^2 + x_1^2 + x_2^2 -8(x_1+x_2))/2},\ \ x=(x_1,x_2)\in[-2,7]^2,
$$
$$
f_4(x)=c_4\e^{-(x_1^2 + (x_2-\alpha(x_1^2-1))^2)/2},\ \ x=(x_1,x_2)\in\mR^2,
$$
where $c_1,c_2,c_3, c_4$ are normalized constants, and
\begin{itemize}
\item $f_1(x)$ represents the 2D-Griewank density function;
\item $f_2(x)$ is the mixture of four normal distributions;
\item  $f_3(x)$ is adopted from \cite[Example 6.4]{RK16};
\item $f_4(x)$ is the Banana-shape distribution with shape parameter $\alpha\in\mR$.
\end{itemize}
We generated 10,000 sample points using the aforementioned algorithm to plot the 2D histogram in Figure \ref{fig:2D1}, \ref{fig:2D2}, \ref{fig:2D3}, \ref{fig:2D4}  below. For the Monte Carlo method, we used 
$M=50$ and 
$n=20,000$ to calculate the integral. The left panel displays the true contour plot, 
the middle panel shows the scatter plot of the sampled points, and the right panel depicts the contour plot of the sampled points.

\begin{figure}[H]
    \centering
    \begin{minipage}[b]{0.3\textwidth}
        \centering
        \includegraphics[width=1.8in, height=1.3in]{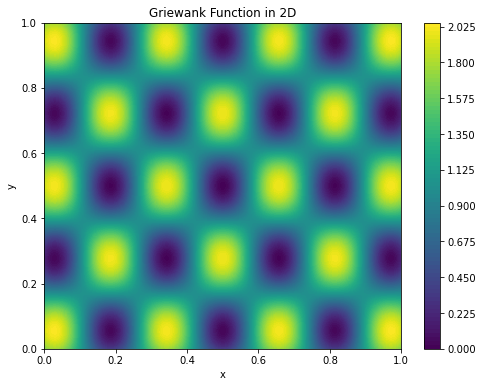}
    \end{minipage}
    \begin{minipage}[b]{0.3\textwidth}
        \centering
        \includegraphics[width=1.8in, height=1.3in]{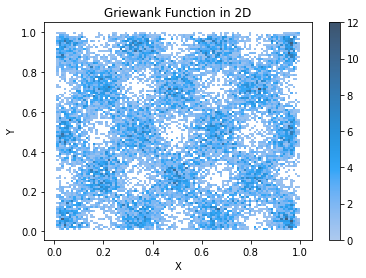}
    \end{minipage}
    \begin{minipage}[b]{0.3\textwidth}
        \centering
        \includegraphics[width=1.8in, height=1.3in]{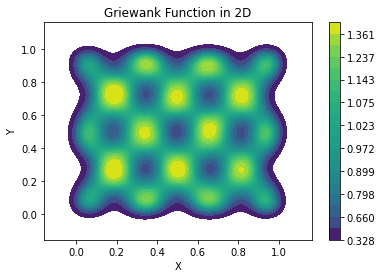}
    \end{minipage}
    \caption{Sampling from $f_1(x)$}
    \label{fig:2D1}
\end{figure}

\begin{figure}[H]
    \centering
    \begin{minipage}[b]{0.3\textwidth}
        \centering
        \includegraphics[width=1.8in, height=1.3in]{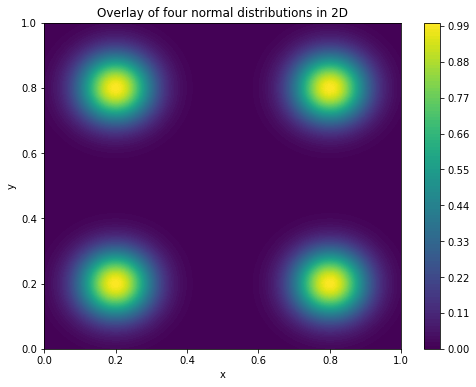}
    \end{minipage}
    \begin{minipage}[b]{0.3\textwidth}
        \centering
        \includegraphics[width=1.8in, height=1.3in]{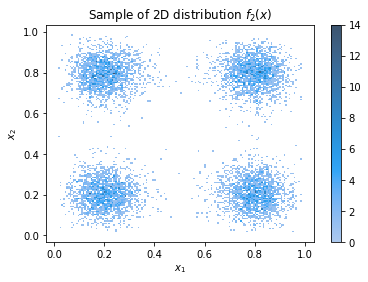}
    \end{minipage}
    \begin{minipage}[b]{0.3\textwidth}
        \centering
        \includegraphics[width=1.8in, height=1.3in]{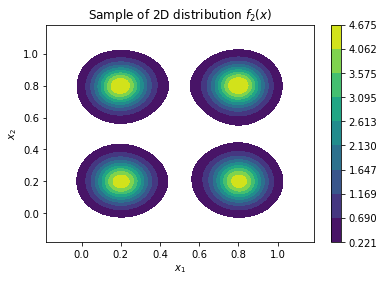}
    \end{minipage}
    \caption{Sampling from $f_2(x)$}
    \label{fig:2D2}
\end{figure}

\begin{figure}[H]
    \centering
    \begin{minipage}[b]{0.3\textwidth}
        \centering
        \includegraphics[width=1.8in, height=1.3in]{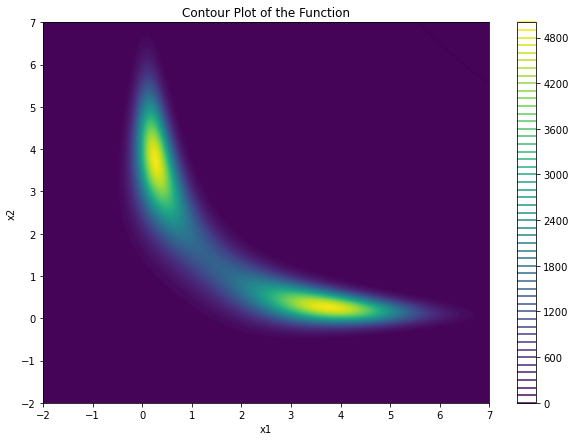}
    \end{minipage}
    \begin{minipage}[b]{0.3\textwidth}
        \centering
        \includegraphics[width=1.8in, height=1.3in]{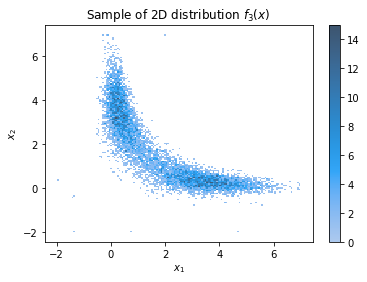}
    \end{minipage}
    \begin{minipage}[b]{0.3\textwidth}
        \centering
        \includegraphics[width=1.8in, height=1.3in]{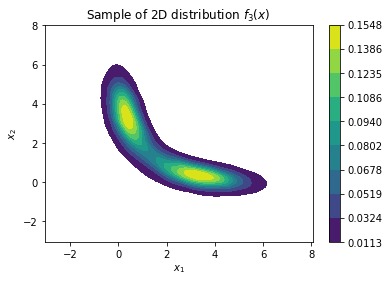}
    \end{minipage}
    \caption{Sampling from $f_3(x)$}
    \label{fig:2D3}
\end{figure}

\begin{figure}[H]
    \centering
    \begin{minipage}[b]{0.3\textwidth}
        \centering
        \includegraphics[width=1.8in, height=1.3in]{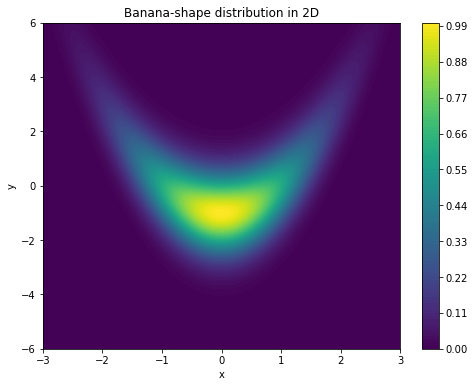}
    \end{minipage}
    \begin{minipage}[b]{0.3\textwidth}
        \centering
        \includegraphics[width=1.8in, height=1.3in]{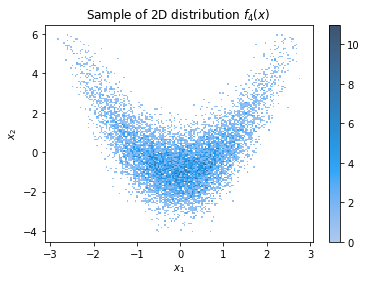}
    \end{minipage}
    \begin{minipage}[b]{0.3\textwidth}
        \centering
        \includegraphics[width=1.8in, height=1.3in]{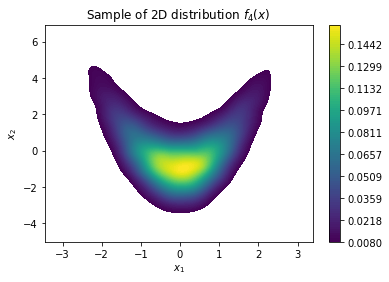}
    \end{minipage}
    \caption{Sampling from $f_4(x)$}
    \label{fig:2D4}
\end{figure}

{\bf Example 3.}  
In this example, we consider the high-dimensional anisotropic funnel distribution. More precisely, we fix a parameter $\alpha > 0$. The anisotropic funnel density is given as follows: for $x=(x_1, x^*_1) \in \mathbb{R} \times \mathbb{R}^{d-1}$,
$$
f(x) = \frac{1}{\sqrt{2\pi}} \exp\left(-\frac{x_1^2}{2}\right) \left(\frac{1}{\sqrt{2\pi \e^{2\alpha x_1}}}\right)^{d-1} \exp\left(-\frac{\|x_1^*\|_2^2}{2\e^{2\alpha x_1}}\right) = \rho_1(x_1) \rho_{\e^{\alpha x_1}}(x_1^*),
$$
where $\rho_\sigma$ is defined in \eqref{DH1}. We generated 10,000 sample points using Algorithm 2 to plot the histogram below when the dimensions are 2, 5, 7, and 10. For the high-dimensional case $d \geq 3$, we project the sampling points onto the first two coordinates. From Figure \ref{fig:2D6} below, one can see that the sampling result deteriorates as the dimension increases.
\begin{figure}[H]
    \centering
    \begin{minipage}[b]{0.3\textwidth}
        \centering
        \includegraphics[width=1.8in, height=1.3in]{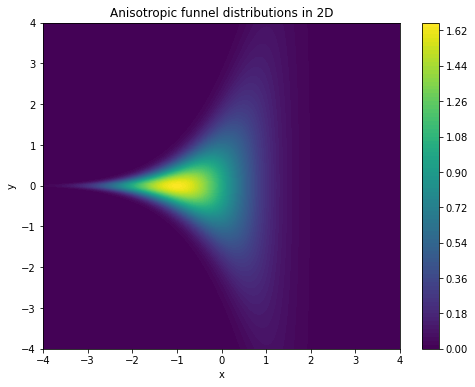}
    \end{minipage}
    \begin{minipage}[b]{0.3\textwidth}
        \centering
        \includegraphics[width=1.8in, height=1.3in]{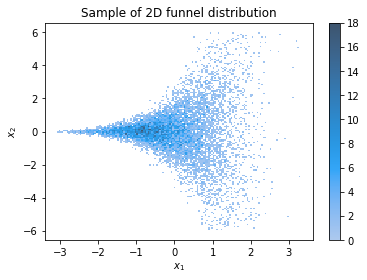}
    \end{minipage}
    \begin{minipage}[b]{0.3\textwidth}
        \centering
        \includegraphics[width=1.8in, height=1.3in]{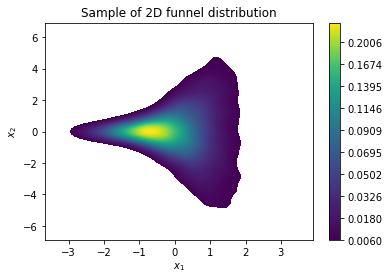}
    \end{minipage}
    \caption{Sampling of  2D-funnel distribution}
    \label{fig:2D5}
\end{figure}

\begin{figure}[H]
    \centering
    \begin{minipage}[b]{0.3\textwidth}
        \centering
        \includegraphics[width=1.8in, height=1.3in]{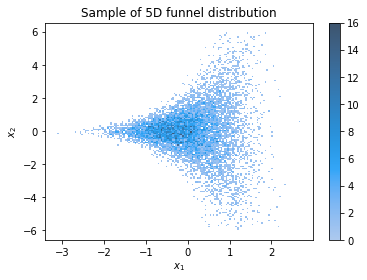}
    \end{minipage}
    \begin{minipage}[b]{0.3\textwidth}
        \centering
        \includegraphics[width=1.8in, height=1.3in]{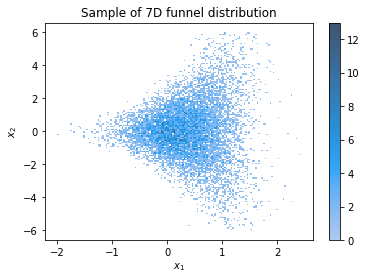}
    \end{minipage}
    \begin{minipage}[b]{0.3\textwidth}
        \centering
        \includegraphics[width=1.8in, height=1.3in]{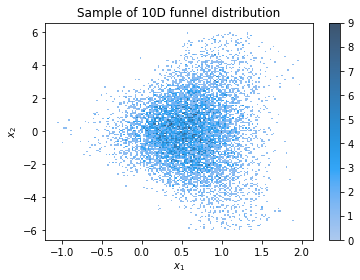}
    \end{minipage}
    \caption{Sampling of 5D, 7D, 10D-funnel distributions}
    \label{fig:2D6}
\end{figure}

Based on the expression of $f$,
one can provide a more explicit expression for $\cD^f_t(x)$ given in \eqref{DX91}.
First of all, by definition and the change of variable we have
\begin{align*}
&\int_{\mR^d}y\rho_{\sigma_t}(x-\beta_ty)f(y)\dif y
=(\beta_t)^{-d-1}\int_{\mR^d}y\rho_{\sigma_t}(x-y)f(y/\beta_t)\dif y\\
&\qquad=\beta_t^{-1}\int_{\mR}\rho_{\sigma_t}(x_1-y_1)\rho_{\beta_t}(y_1)
\left(\int_{\mR^{d-1}}y\rho_{\sigma_t}(x_1^*-y_1^*)\rho_{\beta_t\e^{\alpha y_1}}(y_1^*)\dif y_1^*\right)\dif y_1.
\end{align*}
By elementary calculations, one sees that
$$
\int_{\mR^d}\rho_\sigma(x-y)\rho_{\sigma'}(y)\dif y=\rho_{\sqrt{\sigma^2+\sigma'^2}}(x),
$$
and
$$
\int_{\mR^d}y\rho_\sigma(x-y)\rho_{\sigma'}(y)\dif y=
\frac{x\sigma'^2}{\sigma^2 + \sigma'^2} \rho_{\sqrt{\sigma^2+\sigma'^2}}(x).
$$
Hence,
\begin{align*}
\int_{\mR^{d-1}}y_1\rho_{\sigma_t}(x_1^*-y_1^*)\rho_{\beta_t\e^{\alpha y_1}}(y_1^*)\dif y_1^*
=y_1\rho_{\sqrt{\sigma_t^2+\beta_t^2\e^{2\alpha y_1}}}(x_1^*)
\end{align*} 
and
\begin{align*}
\int_{\mR^{d-1}}y^*_1\rho_{\sigma_t}(x_1^*-y_1^*)\rho_{\beta_t\e^{\alpha y_1}}(y_1^*)\dif y_1^*
=\frac{x^*_1\beta_t^2\e^{2\alpha y_1}}{\sigma_t^2+\beta_t^2\e^{2\alpha y_1}}\rho_{\sqrt{\sigma_t^2+\beta_t^2\e^{2\alpha y_1}}}(x_1^*).
\end{align*} 
Thus
\begin{align*}
\int_{\mR^d}y_1\rho_{\sigma_t}(x-\beta_ty)f(y)\dif y
&=\beta_t^{-1}\int_{\mR}y_1\rho_{\sigma_t}(x_1-y_1)\rho_{\beta_t}(y_1)\rho_{\sqrt{\sigma_t^2+\beta_t^2\e^{2\alpha y_1}}}(x_1^*)\dif y_1\\
&=\int_{\mR}y_1\rho_{\sigma_t}(x_1-\beta_ty_1)\rho_{1}(y_1)\rho_{\sqrt{\sigma_t^2+\beta_t^2\e^{2\alpha\beta_t y_1}}}(x_1^*)\dif y_1\\
&=\mE\Big[ \xi \rho_{\sigma_t}(x_1-\beta_t\xi)\rho_{\sqrt{\sigma_t^2+\beta_t^2\e^{2\alpha\beta_t \xi}}}(x_1^*)\Big]
\end{align*} 
and
\begin{align*}
\int_{\mR^d}y^*_1\rho_{\sigma_t}(x-\beta_ty)f(y)\dif y
&=x^*_1\int_{\mR}\frac{\beta_t\e^{2\alpha y_1}}{\sigma_t^2+\beta_t^2\e^{2\alpha y_1}}\rho_{\sigma_t}(x_1-y_1)\rho_{\beta_t}(y_1)\rho_{\sqrt{\sigma_t^2+\beta_t^2\e^{2\alpha y_1}}}(x_1^*)\dif y_1\\
&=\mE\left[\frac{x^*_1\beta_t\e^{2\alpha \beta_t\xi}}{\sigma_t^2+\beta_t^2\e^{2\alpha\beta_t\xi}}\rho_{\sigma_t}(x_1-\beta_t\xi)\rho_{\sqrt{\sigma_t^2+\beta_t^2\e^{2\alpha \beta_t\xi}}}(x_1^*)\right].
\end{align*} 
Similarly, we have
\begin{align*}
\int_{\mR^d}\rho_{\sigma_t}(x-\beta_ty)f(y)\dif y
&=\int_{\mR}\rho_{\sigma_t}(x_1-y_1)\rho_{\beta_t}(y_1)\rho_{\sqrt{\sigma_t^2+\beta_t^2\e^{2\alpha y_1}}}(x_1^*)\dif y_1\\
&=\mE\left[\rho_{\sigma_t}(x_1-\beta_t\xi)\rho_{\sqrt{\sigma_t^2+\beta_t^2\e^{2\alpha \beta_t\xi}}}(x_1^*)\right].
\end{align*} 
Combining the above calculations, we obtain
\begin{align*}
\cD^{f}_t(x)
&=\frac{\mE\left[\left(\xi, \frac{x^*_1\beta_t\e^{2\alpha \beta_t\xi}}{\sigma_t^2+\beta_t^2\e^{2\alpha\beta_t\xi}}\right)\rho_{\sigma_t}(x_1-\beta_t\xi)\rho_{\sqrt{\sigma_t^2+\beta_t^2\e^{2\alpha \beta_t\xi}}}(x_1^*)\right]}
{\mE\left[\rho_{\sigma_t}(x_1-\beta_t\xi)\rho_{\sqrt{\sigma_t^2+\beta_t^2\e^{2\alpha \beta_t\xi}}}(x_1^*)\right]}.
\end{align*} 
Using this expression, one can devise a similar algorithm for sampling from $f$. Here are the simulated results.
Since this expression only contains the integral with respect to the first variable, it appears to work well for high-dimensional cases.
\begin{figure}[H]
    \centering
    \begin{minipage}[b]{0.3\textwidth}
        \centering
        \includegraphics[width=1.8in, height=1.3in]{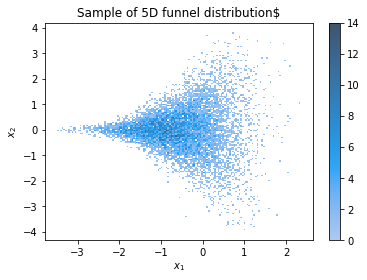}
    \end{minipage}
    \begin{minipage}[b]{0.3\textwidth}
        \centering
        \includegraphics[width=1.8in, height=1.3in]{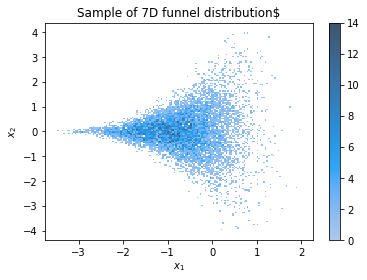}
    \end{minipage}
    \begin{minipage}[b]{0.3\textwidth}
        \centering
        \includegraphics[width=1.8in, height=1.3in]{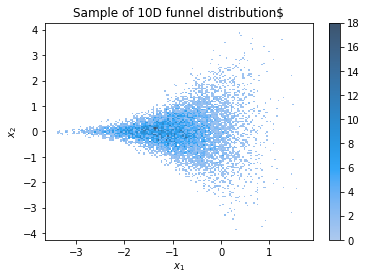}
    \end{minipage}
    \caption{Sampling  of 5D, 7D, 10D-funnel distributions}
    \label{fig:2D7}
\end{figure}

\section{Application to optimizing problems}

Let $U:\mathbb{R}^d \to [0, \infty)$ be a continuous function. In many optimization problems, we need to find the minimum point of $U$, that is,
$$
x_* = \operatorname{argmin}(U(x)).
$$
When $U$ is convex, there are many ways to find its global minimum, for example, the gradient descent method or the stochastic Langevin method. However, when $U$ is non-convex and the dimension $d$ is large, the problem becomes quite challenging. In this case, the probabilistic method would be a good choice. 

Now consider the probability density function
$$
f_\beta(x) = \exp(-\beta U(x)), \quad \beta > 0.
$$
Intuitively, the maximum of $f_\beta$ will correspond to the minimum of $U$, and as $\beta$ becomes large, $f_\beta$ will concentrate on the minimum of $U$. In fact, it has been proved that as $\beta \to \infty$, $f_\beta$ weakly converges to the global minimum of $U$ (see \cite{FW12}).

In this section we use the sampling method in Section \ref{Sec3} to construct an algorithm to seek the minimum of $U$.
We first show the following simple result.
\bt\label{Th41}
Let $f:\mR^d\to[0,\infty)$ be a continuous function.
Suppose that for some $x^*\in\mR^d$,
$$
f^*:=f(x^*)=\sup_{x\in\mR^d}f(x)<\infty.
$$
Let $X_1,\cdots, X_N$ be a sequence of i.i.d. random variables with common distribution $\mu_0(\dif x)=\rho(x)\dif x$. 
Then for any $\eps\in(0,f^*)$ and $N\in\mN$, we have
\begin{align}\label{RA1}
\mE\left|\sup_{n=1,\cdots N}f(X_n)-f^*\right|^2\leq \eps^2+(f^*)^2\e^{-N\delta_\eps},
\end{align}
where $\delta_\eps:=\int_{f(x)>f^*-\eps}\rho(x)\dif x$.
\et
\begin{proof}
Let $Y_n:=f(X_n)$. Then $Y_n$ is a sequence of i.i.d. nonnegative random variables. Let
$$
Y^*_N:=\sup_{n=1,\cdots,N}Y_n
$$
Then we clearly have
$$
\mP(Y^*_N\leq r)=\prod_{n=1}^N\mP(Y_n\leq r)=\mP(Y_1\leq r)^N.
$$
Hence, for any $\eps\in(0,f^*)$,
\begin{align}
\mE\left|Y^*_N-f^*\right|^2
&=2\int^{f^*}_0(f^*-r)\mP(Y^*_N\leq r)\dif r=2\int^{f^*}_0(f^*-r)\mP(Y_1\leq r)^N\dif r\no\\
&=2\left(\int^{f^*}_{f^*-\eps}+\int^{f^*-\eps}_0\right)(f^*-r)\mP(Y_1\leq r)^N\dif r\no\\
&\leq\eps^2 +2\mP(Y_1\leq f^*-\eps)^N\int^{f^*-\eps}_0(f^*-r)\dif r\no\\
&\leq\eps^2+(1-\mP(f(X_1)>f^*-\eps))^N(f^*)^2.\label{SD1}
\end{align}
Note that
$$
\mP(f(X_1)>f^*-\eps)=\int_{f(x)>f^*-\eps}\rho(x)\dif x=\delta_\eps,
$$
and for $x\in(0,1)$,
$$
\log(1-x)\leq -x\Rightarrow (1-x)^N=\e^{N\log(1-x)}\leq \e^{-Nx}.
$$
Substituting these into \eqref{SD1}, we get
$$
\mE\left|Y^*_N-f^*\right|^2\leq\eps^2+(1-\delta_\eps)^N (f^*)^2\leq \eps^2+ (f^*)^2\e^{-N\delta_\eps}.
$$
The proof is complete.
\end{proof}

We have the following corollary.

\bc
In the situation of Theorem \ref{Th41}, if 
$$
K_f:=\frac12\sup_{x\in\mR^d}\sum_{i,j=1,\cdots, d}|\p^2_{ij}f(x)|<\infty,\ \ \lambda:=\int_{\mR^d}f(x)\dif x\in(0,\infty),
$$
then for any $\beta<4/d$ and $N>(f^*)^{2/\beta}$,
$$
\mE\left|\sup_{n=1,\cdots N}f(X_n)-f^*\right|^2\leq N^{-\beta}+(f^*)^2\exp\Big\{-(f^*- N^{-\beta/2}) K_f^{-d/2} N^{1-d\beta/4}/\lambda\Big\},
$$
where $X_n\sim\mu_0(\dif x)=f(x)\dif x/\lambda$. In particular,
$$
\lim_{N\to\infty}\mE\left|\sup_{n=1,\cdots N}f(X_n)-f^*\right|^2=0.
$$
\ec
\begin{proof}
By Taylor's expansion, we have for some $y=\theta x^*+(1-\theta) x$ with $\theta\in[0,1]$,
\begin{align*}
f(x^*)-f(x)&=f(x^*)-f(x)-(x^*-x)\cdot\nabla f(x^*)\\
&=\frac12\tr((x^*-x)\otimes(x^*-x)\cdot\nabla^2 f(y))\\
&\leq\|x^*-x\|^2_\infty K_f.
\end{align*}
Hence,
\begin{align*}
&\int_{f(x)>f^*-\eps}f(x)\dif x
\geq (f^*-\eps){\rm vol}\{x: f(x^*)-f(x)<\eps\}\\
&\geq(f^*-\eps){\rm vol}\{x: \|x^*-x\|^2_\infty K_f<\eps\}=(f^*-\eps)(\eps/K_f)^{d/2}.
\end{align*}
Thus by Theorem \ref{Th41} with $\rho(x)=f(x)/\lambda$, we get
\begin{align*}
\mE\left|\sup_{n=1,\cdots N}f(X_n)-f^*\right|^2\leq \eps^2+(f^*)^2\e^{-N(f^*-\eps) (\eps/K_f)^{d/2}/\lambda}.
\end{align*}
Taking  $\eps=N^{-\beta/2}$, we obtain the desired estimate.
\end{proof}

\br
In the above corollary, the key quantities are $\lambda$ and $K_f$. Smaller $\lambda$ and $K_f$ will lead to faster convergence.
\er

Now we consider the problem of seeking the minimum of a continuous function. 
The following result does not make any regularity assumption about $U$.

\bt
Let $U:\mR^d\to[0,\infty)$ be a continuous function. Suppose that  there is a  minimum point $x_*\in\mR^d$ so that 
$$
U_*:=\inf_{x\in\mR^d}U(x)=U(x_*)=0,
$$
and for some $0<\kappa_0\leq\kappa_1<\infty$ and $\ell\geq 0$,
\begin{align}\label{SP1}
\kappa_0\|x-x_*\|^2_2\1_{\kappa_0\|x-x_*\|^2_2\geq\ell}\leq U(x)\leq\kappa_1\|x-x_*\|^2_2.
\end{align}
Fix $\beta\geq 1$ and
let $X_1,\cdots, X_N$ be a sequence of i.i.d. random variables with common distribution 
$\mu_0(\dif x)\varpropto\e^{-\beta U(x)}\dif x$.
Then for any $\eps>0$, it holds that
\begin{align}\label{RA03}
\mE\left(\e^{-\inf_{n=1,\cdots N}U(X_n)}-1\right)^2\leq \eps^2+\e^{-N\delta_{\beta,\eps}(\ell)},
\end{align}
where for $\gamma(d,a):=\int^{a}_0 u^{d-1}\e^{-u}\dif u$ (lower incomplete gamma function),
$$
\delta_{\beta,\eps}(\ell):=\frac{(\kappa_0/\kappa_1)^{d/2}\gamma(\tfrac{d}{2},\beta\eps)}
{\Gamma(\tfrac d2)-\gamma(\tfrac d2,\ell\beta)+(\ell\beta)^{d/2}/d}.
$$
\et
\begin{proof}
Applying  \eqref{RA1} to $f(x)=\e^{-U(x)}$ and $\rho(x)=\e^{-\beta U(x)}/\int_{\mR^d}\e^{-\beta U(x)}\dif x$, we have
\begin{align}\label{RA01}
\mE\left(\e^{-\inf_{n=1,\cdots N}U(X_n)}-1\right)^2=
\mE\left|\sup_{n=1,\cdots N}\e^{-U(X_n)}-1\right|^2\leq \eps^2+\e^{-N\delta_{\beta,\eps}},
\end{align}
where 
$$
\delta_{\beta,\eps}=\int_{\e^{-U(x)}>1-\eps}\e^{-\beta U(x)}\dif x/\int_{\mR^d}\e^{-\beta U(x)}\dif x.
$$
Let us first estimate the numerator of $\delta_{\beta,\eps}$. 
Noting that
$$
1-\e^{-U(x)}\leq U(x)\leq \kappa_1\|x-x_*\|^2_2,
$$
we have
$$
\Big\{x: \kappa_1\|x-x_*\|^2_2<\eps\Big\}\subset\Big\{x: U(x)<\eps\Big\}\subset\Big\{x: \e^{-U(x)}>1-\eps\Big\}.
$$
Hence, by \eqref{SP1} and the change of variable,
\begin{align*}
\int_{\e^{-U(x)}>1-\eps}\e^{-\beta U(x)}\dif x&\geq\int_{\kappa_1\|x\|^2_2<\eps}\e^{-\beta U(x)}\dif x
\leq\int_{\kappa_1\|x-x_*\|^2_2<\eps}\e^{-\beta\kappa_1\|x-x_*\|^2_2}\dif x\\
&=\int_{\|x\|^2_2<\beta\eps}\frac{\e^{-\|x\|^2_2}}{(\beta\kappa_1)^{\frac d2}}\dif x
=\frac{\pi^{\frac d2}\gamma(\tfrac{d}{2}, \beta\eps)}{(\beta\kappa_1)^{\frac d2}\Gamma(\frac d2)},
\end{align*}
where we have used that for $a>0$, 
\begin{align}\label{DW1}
\int_{\|x\|_2<a}\e^{-\|x\|^2_2}\dif x=\frac{2\pi^{\frac d2}}{\Gamma(\frac d2)}\int^a_0r^{d-1}\e^{-r^2}\dif r
=\frac{\pi^{\frac d2}\gamma(\tfrac{d}{2}, a^2)}{\Gamma(\frac d2)}.
\end{align}
Next, we look at the denominator of $\delta_{\beta,\eps}$. We decompose it into two parts:
$$
\int_{\mR^d}\e^{-\beta U(x)}\dif x=\int_{\kappa_0\|x-x_*\|^2_2>\ell}\e^{-\beta U(x)}\dif x
+\int_{\kappa_0\|x-x_*\|^2_2\leq\ell}\e^{-\beta U(x)}\dif x=:I_1(\ell)+I_2(\ell).
$$
For $I_1(\ell)$, by \eqref{SP1}, the change of variable and \eqref{DW1}, we have
\begin{align*}
I_1(\ell)\leq\int_{\kappa_0\|x-x_*\|^2_2>\ell}\e^{-\beta\kappa_0\|x-x_*\|^2_2}\dif x
&=\int_{\kappa_0\|x\|^2_2>\ell\beta}\frac{\e^{-\|x\|^2_2}}{(\beta\kappa_0)^{\frac d2}}\dif x
=\frac{\pi^{\frac d2}}{(\beta\kappa_0)^{\frac d2}}\left(1-\frac{\gamma(\tfrac d2,\ell\beta)}{\Gamma(\tfrac d2)}\right).
\end{align*}
For $I_2(\ell)$, we clearly have
$$
I_2(\ell)\leq {\rm vol}(\{x: \kappa_0\|x-x_*\|^2_2<\ell\})=\tfrac{2\pi^{d/2}}{\Gamma(d/2+1)}(\ell/\kappa_0)^{d/2}=\tfrac{\pi^{d/2}}{d\Gamma(d/2)}(\ell/\kappa_0)^{d/2}.
$$
Hence,
$$
\int_{\mR^d}\e^{-\beta U(x)}\dif x\leq
\frac{\pi^{\frac d2}}{(\beta\kappa_0)^{\frac d2}}\left(1-\frac{\gamma(\tfrac d2,\ell\beta)}{\Gamma(\tfrac d2)}
+\frac{(\ell\beta)^{d/2}}{d\Gamma(\frac d2)}\right).
$$
Combining the above calculations, we obtain
$$
\delta_{\beta,\eps}\geq\left(\frac{\kappa_0}{\kappa_1}\right)^{\frac d2}
\frac{\gamma(\tfrac{d}{2}, \beta\eps)}{\Gamma(\frac d2)-\gamma(\tfrac d2,\ell\beta)+(\ell\beta)^{d/2}/d}.
$$
The proof is complete.
\end{proof}

\br
By the definition of $\delta_\beta(\ell)$, it is easy to see that
$$
\delta_{\beta,\eps}(0)=(\kappa_0/\kappa_1)^{\frac d2}\gamma(\tfrac{d}{2},\beta\eps)/\Gamma(\tfrac d2)\stackrel{\beta\to\infty}{\to} (\kappa_0/\kappa_1)^{\frac d2},
$$
and for fixed $\eps,\ell>0$,
$$
\lim_{\beta\to\infty}\delta_{\beta,\eps}(\ell)=0.
$$
Here is the picture of the domain $\{(x,y): \|x\|_2^2 1_{\|x\|_2>1.5} \leq y \leq 2\|x\|_2^2\}$ where the function $U$ stays.
\begin{center}
\begin{tikzpicture}
    \begin{axis}[
        domain=-2:2, 
        samples=500,
        xlabel={$x$},
        ylabel={$y$},
        axis lines=middle,
        legend pos=outer north east,
        grid=both,
        ymin=0, ymax=8,
        xmin=-2, xmax=2,
        width=12cm,
        height=6cm
    ]

    \addplot [
        name path=A,
        domain=-2:-1.5,
        samples=50,
        thick,
        blue,
        dashed
    ] {x^2};

    \addplot [
        name path=B,
        domain=1.5:2,
        samples=50,
        thick,
        blue,
        dashed
    ] {x^2};

    \addplot [
        name path=C,
        domain=-2:2,
        samples=500,
        thick,
        red
    ] {2*x^2};

    \addplot [
        name path=D,
        domain=-1.5:1.5,
        samples=500,
        thick,
        red
    ] {2*x^2};

    \addplot [
        name path=E,
        domain=-1.5:1.5,
        samples=500,
        draw=none
    ] {0};

    \addplot [
        fill=gray, 
        fill opacity=0.5
    ] fill between[
        of=A and C,
        soft clip={domain=-2:-1.5}
    ];

    \addplot [
        fill=gray, 
        fill opacity=0.5
    ] fill between[
        of=B and C,
        soft clip={domain=1.5:2}
    ];

    \addplot [
        fill=gray, 
        fill opacity=0.5
    ] fill between[
        of=E and D,
        soft clip={domain=-1.5:1.5}
    ];

    \end{axis}
\end{tikzpicture}
\end{center}
In particular,  when $\ell=0$, condition \eqref{SP1} means that $U$ lies in the middle of two parabolic line.
\er

Based on the above theorem, we devise a new algorithm to find the minimum point of a function $U$.
More precisely, we fix the number of search iterations $M$ and the number of sample points $N$. We initialize the intensity $\beta \geq 1$, the scale parameter $\alpha \geq 1$, and the initial position $x_*$ of the minimum point of $U$. Using Algorithm 1, we sample $N$ points $\mathcal{S}_N$ from the distribution $\mu_0(\mathrm{d} x) \propto \exp(-\beta U(\alpha(x - x_*) + x_*))$ within the domain $\{x: \|x - x_*\|_\infty \leq 1\}$. Next, we calculate the minimum value of $U$ over $\mathcal{S}_N$ and update $x_*$ to the corresponding minimum point of $U$ over $\mathcal{S}_N$. Finally, by appropriately increasing $\beta$ and decreasing $\alpha$, we repeat the above procedure $M$ times to search for the minimum point.
Here is the concrete algorithm.

\noindent\rule{\linewidth}{1 pt}

{\bf Algorithm 3.} Find the minimum of a function through sampling the density function

\noindent\rule{\linewidth}{0.8 pt}

\begin{enumerate}[1.]
\item {\bf Data}: Give a nonnegative function $U$, sample points number $N$, search times $M$.
\item {\bf Initial data:} Intensity $\beta=1$, scale $\alpha=10$ and the minimum point $x_*=0$.
\item {\bf For} $j=0$ {\bf to} $M-1$
\begin{itemize}
\item $U_\alpha(x)=U(\alpha(x-x_*)+x_*)$.
\item Use Algorithm 2 to generate $N$-sample points $\cS_N$ from $\mu_0(\dif x)\varpropto\e^{-\beta U_\alpha(x)}$ in the domain 
$$
\{x: \|x-x_*\|_\infty\leq 1\}.
$$
\item Update $x_*\leftarrow {\rm argmin}_{\cS_N}U_\alpha$, $\beta\leftarrow \beta*(5+j)$, $\alpha\leftarrow\alpha/\ln(3+j)$.
\end{itemize}
\item {\bf Output}: $x_*$ and $U(x_*)$.
\end{enumerate}
\noindent\rule{\linewidth}{1 pt}

In the above algorithm, with each iteration, we increase $\beta$ to $(5 + j)\beta$ and decrease the scale $\alpha$ to $\alpha / \log(3 + j)$, where $j$ is the current iteration number. The factors $5 + j$ and $\log(3 + j)$ can be adjusted for different tasks.

In the following examples, $U_1$ is the Griewank function, $U_2$ is the Rosenbrock function, $U_3$ is the Ackley function, $U_4$ is the Rastrigin function, $U_5$ is a quadratic function, and $U_6$ is the sum of two Gaussian functions. The functions $U_1$, $U_2$, $U_3$, and $U_4$ are commonly used as test functions in various optimization algorithms.
We apply {\bf Algorithm 3} to these functions with the parameters set as follows: the number of sample points $N = 10$, 
the number of search times $M = 5$, and in {\bf Algorithm 2}, we set $n = 50000$ and the iteration time to $30$ for generating samples.
We would like to highlight that we sample only 10 points to determine the minimum points.

$$
U_1(x)=\tfrac{x_1^2 + x_2^2}{ 4000}-\cos(x_1)\cos(\tfrac{x_2}{\sqrt 2})+1,\qquad\qquad
U_2(x)=(1-x_1)^2+100(x_1^2-x_2)^2
$$
    \begin{minipage}{0.5\textwidth}
        \centering
        \begin{tabular}{|c|c|c|c|}
        \hline
         $M$ & $x_*$=argmin($U_1$) & $\min(U_1)$\\
        \hline
        1 & 0.355809, 0.307940 & 0.125017362963\\
        \hline
        2& 0.355809, 0.307940 & 0.125017362963\\
        \hline
        3&0.355809, 0.307940 & 0.125017362963\\
        \hline
        4& 0.499977, 0.499881& 1.518449525e-06\\
        \hline
        5&0.500011, 0.500081& 7.158456633e-07 \\
        \hline
        \end{tabular}
    \end{minipage}
    \hfill
    \begin{minipage}{0.5\textwidth}
        \centering
        \begin{tabular}{|c|c|c|c|}
        \hline
         $M$ & $x_*$=argmin($U_2$) & $\min(U_2)$\\
        \hline
        1 & 1.066735,  1.156329 & 0.038330061783\\
        \hline
        2&1.005933, 1.019098 & 0.005213454835\\
        \hline
        3&1.009863,  1.021328&0.000323506548\\
        \hline
        4& 0.999039, 0.998426& 1.296778123e-05\\
        \hline
        5&1.000203,  1.000409& 4.184867527e-08 \\
        \hline
        \end{tabular}
    \end{minipage}
$$
U_3(x)=20-\e^{-\frac{\|x\|_2}{5\sqrt 2}}+\e-\e^{\sum_{i=1,2}\frac{\cos(2\pi x_i)}2},\qquad\qquad
U_4(x)=20+\frac{\|x\|^2_2}2-10\sum_{i=1,2}\cos(2\pi x_i).
$$
\begin{minipage}{0.5\textwidth}
        \centering
        \begin{tabular}{|c|c|c|c|}
        \hline
         $M$ & $x_*$=argmin($U_4$) & $\min(U_4)$\\
        \hline
        1 & -0.001467, -0.002903& 0.00948314937\\
        \hline
        2&-0.000869,  -0.000305& 0.00263084907\\
        \hline
        3&2.209e-05, -1.850e-04 & 0.00052806219\\
        \hline
        4& -3.708e-06,  8.937e-06& 2.7370639e-05\\
        \hline
        5&-3.708e-06,  8.937e-06& 2.7370639e-05 \\
        \hline
        \end{tabular}
    \end{minipage}
    \hfill
    \begin{minipage}{0.5\textwidth}
        \centering
        \begin{tabular}{|c|c|c|c|}
        \hline
         $M$ & $x_*$=argmin($U_2$) & $\min(U_2)$\\
        \hline
        1 & -0.000916, -0.001308 & 0.0005064938\\
        \hline
        2& -0.000148, -0.000496 & 5.3346907e-05\\
        \hline
        3&4.122e-05, -2.363e-04 & 1.14172916e-05\\
        \hline
        4& -9.402e-05,  6.881e-05 & 2.69342732e-06\\
        \hline
        5&1.507e-06, 1.381e-05 & 3.83101905e-08\\
        \hline
        \end{tabular}
    \end{minipage}
$$
U_5(x)=\|x-0.3\|_2^2+\|x-0.1\|_2^2,\qquad\qquad\qquad U_6(x)=2-\e^{-\|x-0.3\|_2^2}-\e^{-\|x+0.3\|_2^2}.
$$
\begin{minipage}{0.5\textwidth}
        \centering
        \begin{tabular}{|c|c|c|c|}
        \hline
         $M$ & $x_*$=argmin($U_4$) & $\min(U_4)$\\
        \hline
        1 & 0.181246, 0.173710 & 0.042085701865\\
        \hline
        2& 0.210050, 0.200877 & 0.040203557243\\
        \hline
        3&0.203378, 0.200219 & 0.040022920593\\
        \hline
        4& 0.199966, 0.200430 & 0.040000373147\\
        \hline
        5&0.199966, 0.200430 & 0.040000373147 \\
        \hline
        \end{tabular}
    \end{minipage}
    \hfill
    \begin{minipage}{0.5\textwidth}
        \centering
        \begin{tabular}{|c|c|c|c|}
        \hline
         $M$ & $x_*$=argmin($U_2$) & $\min(U_2)$\\
        \hline
        1 & -0.040019,  -0.047588 & 0.33360633586\\
        \hline
        2&0.015872, 0.004686 & 0.32979003108\\
        \hline
        3&0.005617,  0.001760 & 0.32950110148\\
        \hline
        4& 0.005617,  0.001760 & 0.32950110148\\
        \hline
        5&0.005617,  0.001760 & 0.32950110148 \\
        \hline
        \end{tabular}
    \end{minipage}

\vspace{5mm}

{\bf Acknowledgement:} 
I would like to express my deep thanks to Zimo Hao, Zhenyao Sun and Rongchan Zhu for their quite useful mathematical discussions, and to Mingyang Lai, Qi Meng and Dianpeng Wang  for their exceptional assistance with coding.


\begin{thebibliography}{99}

\bibitem{ABE23}M. S. Albergo, N. M. Boffi, and E. Vanden-Eijnden. 
Stochastic Interpolants: A Unifying Framework for Flows and Diffusions. 
 \textit{arXiv preprint arXiv:2303.08797v3}, 2023.

\bibitem{And82} B. Anderson. Reverse-time diffusion equation models. \textit{Stochastic Process. Appl.}, 12(3):313–326, May 1982.

\bibitem{B06} C. M. Bishop. \textit{Pattern Recognition and Machine Learning}. Information Science and Statistics, Springer, 2006.

\bibitem{CCLLZ22} S. Chen, S. Chewi, J. Li, Y. Li, A. Salim, and A. R. Zhang. Sampling is as easy as learning the score: Theory for diffusion models with minimal data assumptions. \textit{arXiv preprint arXiv:2209.11215}, 2022.

\bibitem{CLL23} H. Chen, H. Lee, and J. Lu. Improved analysis of score-based generative modeling: User-friendly bounds under minimal smoothness assumptions. In \textit{International Conference on Machine Learning}, pages 4735–4763. PMLR, 2023.

\bibitem{CMFW24}M. Chen, S. Mei, J. Fan and M. Wang.
An Overview of Diffusion Models: Applications, Guided Generation, Statistical Rates and Optimization. arXiv:2404.07771v1.

\bibitem{CLTX23} X. Cheng, J. Lu, Y. Tan, and Y. Xie. Convergence of flow-based generative models via proximal gradient descent in Wasserstein space. \textit{arXiv preprint arXiv:2310.17582v1}, 2023.

\bibitem{DWMG15} J. Sohl-Dickstein, E. Weiss, N. Maheswaranathan, and S. Ganguli. Deep unsupervised learning using nonequilibrium thermodynamics. In \textit{International Conference on Machine Learning}. PMLR, 2015, pp. 2256–2265.

\bibitem{ELS20} R. Eldan, J. Lehec, and Y. Shenfeld. Stability of the logarithmic Sobolev inequality via the F\"ollmer process. \textit{Annales de l’Institut Henri Poincar\'e, Probabilit\'es et Statistiques}, 56:2253–2269, 2020.

\bibitem{FW12} M. I. Freidlin and A. D. Wentzell. \textit{Random Perturbations of Dynamical Systems}. Third Edition. Translated by Joseph Szücs, Springer, 2012.

\bibitem{F84} H. F\"ollmer. Time reversal on Wiener space. In \textit{Stochastic Processes—Mathematics and Physics}, Lecture Notes in Math., 1158, pages 119-129, Springer, Berlin, 1986.

\bibitem{FG15} N. Fournier and A. Guillin. On the rate of convergence in Wasserstein distance of the empirical measure. \textit{Probab. Theory Relat. Fields}, 162:707-738, 2015.

\bibitem{GHJ23}Y. Gao, J. Huang and Y. Jiao. Gaussian Interpolation Flows. \textit{arXiv preprint arXiv:2311.11475v1}.

\bibitem{GCSR04} A. Gelman, J. B. Carlin, H. S. Stern, and D. B. Rubin. \textit{Bayesian Data Analysis}. Chapman \& Hall/CRC Press, London, New York, 2004.

\bibitem{GMGD24} L. Grenioux, N. Maxence, M. Gabri\'e, and A. Durmus. Stochastic Localization via Iterative Posterior Sampling. \textit{arXiv preprint arXiv:2402.10758v2}, 2024.

\bibitem{HP86} U. G. Haussmann and E. Pardoux. Time reversal of diffusions. \textit{The Annals of Probability}, 14(4):1188-1205, 1986.

\bibitem{HTF08} T. Hastie, R. Tibshirani, and J. Friedman. \textit{The Elements of Statistical Learning: Data Mining, Inference, and Prediction}. Springer Series in Statistics, 2008.

\bibitem{HJA20} J. Ho, A. Jain, and P. Abbeel. Denoising Diffusion Probabilistic Models. In \textit{Proceedings of the 37th International Conference on Machine Learning (ICML)}, 2020.

\bibitem{HJK21} J. Huang, Y. Jiao, L. Kang, X. Liao, J. Liu, and Y. Liu. Schr\"odinger-F\"ollmer sampler: Sampling without ergodicity. \textit{arXiv preprint arXiv:2106.10880}, 2021.

\bibitem{HDMZ24} X. Huang, H. Dong, Y. Hao, Y. Ma, and T. Zhang. Reverse Diffusion Monte Carlo. In \textit{The Twelfth International Conference on Learning Representations}, 2024.

\bibitem{KW08} M. H. Kalos and P. A. Whitlock. \textit{Monte-Carlo Method}. WILEY-VCH Verlag GmbH \& Co., KGaA, Weinheim, 2008.

\bibitem{L10} S. L. Lohr. \textit{Sampling: Design and Analysis}, 2nd ed. Advanced Series, Boston, 2010.

\bibitem{LLT23} H. Lee, J. Lu, and Y. Tan. Convergence of score-based generative modeling for general data distributions. In \textit{International Conference on Algorithmic Learning Theory}, pages 946–985. PMLR, 2023.

\bibitem{LZ22} C. Lu, Y. Zhou, F. Bao, J. Chen, C. Li, and J. Zhu. DPM-Solver: A Fast ODE Solver for Diffusion Probabilistic Model Sampling in Around 10 Steps. In \textit{36th Conference on Neural Information Processing Systems (NeurIPS 2022)}.

\bibitem{M12} K. P. Murphy. \textit{Machine Learning: A Probabilistic Perspective}. The MIT Press, Cambridge, Massachusetts, London, England, 2012.

\bibitem{RK16} Y. R. Rubinstein and D. P. Kroese. \textit{Simulation and the Monte Carlo Method}. Wiley Series in Probability and Statistics, 2016.

\bibitem{RXZ20} M. R\"ockner, L. Xie, and X. Zhang. Superposition principle for non-local Fokker–Planck–Kolmogorov operators. \textit{Probability Theory and Related Fields}, 178:699-733, 2020.

\bibitem{SD21} Y. Song, J. Sohl-Dickstein, D. P. Kingma, A. Kumar, S. Ermon, and B. Poole. Score-based generative modeling through stochastic differential equations. In \textit{International Conference on Learning Representations}, 2021.

\bibitem{SE19} Y. Song and S. Ermon. Generative modeling by estimating gradients of the data distribution. In \textit{Advances in Neural Information Processing Systems}, pages 11895–11907, 2019.

\bibitem{TL23} L. Triplett and J. Lu. Diffusion Methods for Generating Transition Paths. \textit{arXiv preprint arXiv:2309.10276v1}, 2023.

\bibitem{YZS23} L. Yang, Z. Zhang, Y. Song, S. Hong, and R. Xu et al. Diffusion Models: A Comprehensive Survey of Methods and Applications. \textit{ACM Computing Surveys}, 56(4): Article 105, November 2023.

\bibitem{Zh23} X. Zhang. Compound Poisson particle approximation for McKean-Vlasov SDEs. \textit{arXiv preprint arXiv:2306.06816}, 2023.

\end{thebibliography}
\end{document}